%% file: acs26_arxiv.tex
\pdfoutput=1
\documentclass[11pt,a4paper]{article}
\usepackage[T1]{fontenc}
\usepackage{lmodern}
\usepackage{graphicx}
\usepackage{amsmath,amssymb}
\usepackage{booktabs}
\usepackage{multirow}
\usepackage{url}
\usepackage{array}
\usepackage[margin=25mm]{geometry}
\usepackage{verbatim}
\usepackage[ruled,vlined,linesnumbered]{algorithm2e}
\usepackage{tikz}
\usetikzlibrary{positioning,arrows.meta,fit,backgrounds}
\newcommand{\ee}[2]{$#1\times10^{#2}$}
\newcommand{\spd}[1]{$\times#1$}
\newcommand{\TS}{\textsc{TwoSum}}
\newcommand{\FTS}{\textsc{FastTwoSum}}
\newcommand{\doi}[1]{\texttt{doi:#1}}
\newcommand{\TP}{\textsc{TwoProd}}
\SetKwFunction{TwoSum}{TwoSum}\SetKwFunction{FastTwoSum}{FastTwoSum}\SetKwFunction{TwoProd}{TwoProd}
\SetKwInOut{Input}{input}\SetKwInOut{Output}{output}
\IncMargin{1em}
\usepackage[hidelinks]{hyperref}

\title{Performance evaluation of branch-free fused multiply-add algorithms\\
for multicomponent multiple-precision floating-point arithmetic}
\author{Tomonori Kouya\\
\small Faculty of Science and Engineering, Otemon Gakuin University}
\date{}

\begin{document}
\maketitle

\begin{abstract}
Multicomponent multiple-precision arithmetic is constructed from existing floating-point operations using error-free transformations (EFTs). Its performance can be improved by employing branch-free algorithms that eliminate conditional branches. Zhang and Aiken proposed branch-free addition and multiplication algorithms for double-word (DW), triple-word (TW), and quad-word (QW) arithmetic. We implemented their TW and QW algorithms, for which substantial performance improvements over conventional algorithms were expected, and demonstrated their effectiveness. In this paper, we propose branch-free fused multiply-add (FMA) algorithms for DW, TW, and QW arithmetic. The proposed algorithms integrate multiplication and addition into a single computational network and require fewer arithmetic operations than separately performing branch-free multiplication and addition. The anchor-relative error bounds, the preconditions of all FastTwoSum operations, and the non-overlapping properties of the outputs are mechanically verified using FPANVerifier, and input-relative error bounds are subsequently derived analytically. Benchmark results on CPUs and GPUs show that the proposed algorithms provide performance improvements in many compute-intensive cases, including division, square root, and basic linear algebra kernels, while maintaining accuracy comparable to that of the existing branch-free algorithms.
\end{abstract}


%--------------------------------------%
% Introduction
%--------------------------------------%
\section{Introduction}\label{sec:intro}

One approach to ensuring the accuracy of computed results in finite-precision
floating-point arithmetic is increasing the significand length.
Among these approaches, multicomponent multiple-precision arithmetic
(double-word (DW), triple-word (TW), quad-word (QW); denoted in our implementation as DD/TD/QD and,
for a single-precision base, DS/TS/QS), represents numbers as the unevaluated sum of several binary64 or binary32
components and handles the carries exactly by error-free transformations
(EFTs). This approach generally achieves higher performance than arbitrary-precision libraries
because it directly exploits existing floating-point units.
However, its interior, in particular the normalization stage, contains if-branches
that depend on the magnitudes of the values, which are the main causes of
performance loss due to lane or warp divergence when the code is vectorized
for SIMD or ported to a GPU.

Zhang and Aiken proposed branch-free algorithms that eliminate these branches for
addition and multiplication in DW/TW/QW, and demonstrated high performance together
with a machine proof of the error bounds using an SMT solver\cite{ZhangAiken25SC,ZhangAiken25CAV}.
Among these, we have implemented TW and QW, for which substantial performance gains over
the existing algorithms can be expected, and have confirmed that they contribute
to acceleration\cite{tkouya_bf}.

To further improve performance, this study proposes a fused multiply-add
(FMA) $z\leftarrow x\cdot y+c$ for DW/TW/QW, which we constructed using
Claude Code Opus 4.8.
A hardware FMA is commonly used to implement multiply--accumulate (MAC) operations;
most current CPUs and GPUs provide it as a single instruction, and it is widely used
in basic linear algebra kernels.
In this paper, however, the term FMA denotes a multicomponent multiply--add operation
that integrates multiplication and addition into a single computational network rather
than executing them as two independent multicomponent operations.  Thus, the proposed
operation should not be confused with the IEEE~754 single-rounding hardware FMA used
internally as an EFT primitive.

The proposed FMA needs fewer operations than a combination of a branch-free addition
and multiplication (DW $17$ / TW $72$ / QW $176$ flops against $29/96/209$ flops),
its anchor-relative error bound and the validity of every FastTwoSum premise have been
machine-proved with FPANVerifier, an input-relative error bound has been derived by an
analytic transformation based on the non-overlapping condition of the inputs, and it is
bitwise commutative with respect to the exchange $x\leftrightarrow y$.
As a reference that measures the attainable accuracy, we also provide an Exact FMA with a
binomial-type error bound in which the result-relative term dominates within a certain
cancellation range.
For the Exact FMA we show that with $n=K$ distillation sweeps
an input-relative residual term remains of the same order as the result-relative
term, so that under deep cancellation the name ``Exact'' is not justified; by unifying
the rule to $n=K{+}1$ we confirm that the measured maximum error of a single operation
becomes at most $0.5\,\varepsilon_K$ (an error level corresponding to correct rounding
within the tested range) for all six types.

The remainder of this paper is organized as follows.
Section~\ref{sec:algo} describes the structures of the proposed algorithms and their
machine-proved error bounds, and Section~\ref{sec:related} contrasts them with existing
work on FMA. Section~\ref{sec:bench} presents the benchmark results for the two
environments, GB10 (Arm) and H100 (x86)---the acceleration of division and square root,
and the basic linear algebra kernels AXPY/GEMV/GEMM (serial, OpenMP and GPU parallel).
Section~\ref{sec:concl} presents conclusions and future work directions.

%--------------------------------------%
% The proposed FMA
%--------------------------------------%
\section{Branch-free FMA algorithms for DW, TW, and QW}\label{sec:algo}

The proposed branch-free FMA is a fused operation that computes
$z\leftarrow x\cdot y+c$ in a single stage for $K$-component numbers
$x = (x_0, x_1, \ldots, x_{K-1})$ $(=\sum_{i=0}^{K-1} x_i)$,
$y = (y_0, y_1, \ldots, y_{K-1})$ and $c = (c_0, c_1, \ldots, c_{K-1})$
($K=2,3,4$ for DW/TW/QW).  It has the following properties.
\begin{itemize}
\item \textbf{Commutativity}: the operation is commutative in the two multiplicative
      operands $x$ and $y$; $\mathrm{fma}(x,y,c)$ and $\mathrm{fma}(y,x,c)$ agree
      bitwise for all six types
      (200{,}000 trials each, no mismatch).
\item \textbf{Machine-proved error bound}: FPANVerifier mechanically verifies the
      anchor-relative error bounds and the validity of every \FTS{} premise for all
      precisions $p\ge8$ (including binary32/64/128), where $p$ denotes the significand
      length of the underlying binary floating-point format; an analytic transformation
      based on the non-overlapping conditions of the inputs $x$, $y$, and $c$ then yields
      the input-relative bounds $35u^2$, $187u^3$, and $822u^4$, respectively, times
      $(|xy|+|c|)$
      (Section~\ref{sec:verify}).
\item \textbf{Reduced operation count}: it computes with $17/72/176$ flops, fewer than
      the $29/96/209$ flops of a branch-free multiplication combined with an addition
      (Section~\ref{sec:flops}).
\item \textbf{Error semantics}: the proposed version has an error bound relative to
      $(|xy|+|c|)$, whereas the Exact version has a binomial-type error bound consisting
      of a result-relative term and an input-relative residual term, in which the
      result-relative term dominates within a certain cancellation range
      (Section~\ref{sec:exact}).
      In this study, the distillation sweeps of the Exact version are unified to
      $n=K{+}1$, such that the result-relative property is preserved down to deep
      cancellation ($38$--$44$ bits at $p{=}53$).
\end{itemize}

The proposed FMA was built using only three error-free transformations (EFTs)
(Table~\ref{tab:eft}).

These EFTs represent a sum or product exactly as the sum of two floating-point numbers.
\TS{} goes back to Knuth\cite{Knuth}, while \FTS{} and the original \TP{} construction
go back to Dekker\cite{Dekker71}.  Here we use the FMA-based two-operation form of \TP{}.
\FTS{} is three flops cheaper than \TS{}, but it is exact only when
$\mathrm{exponent}(a)\ge\mathrm{exponent}(b)$; this premise is the crux of the machine
proofs in this paper.
\TP{} obtains the rounding error of a product in $2$ flops by using the hardware FMA.

\begin{table}[h]\centering
\caption{Building blocks and their flop counts (a hardware FMA counts as $1$)}\label{tab:eft}
{\small\begin{tabular}{ll}
\toprule
operation & flops \\
\midrule
$(s,e):=\TwoSum(a,b)$      & 6 \\
\hspace{2em}$s := \mathrm{fl}(a + b)$;\ $v := \mathrm{fl}(s - a)$ & \\
\hspace{2em}$e := \mathrm{fl}(\mathrm{fl}(a - \mathrm{fl}(s - v)) + \mathrm{fl}(b - v))$ \\
\midrule
$(s,e):=\FastTwoSum(a,b)$  & 3 \\
\hspace{2em}$s := \mathrm{fl}(a + b)$ &\\
\hspace{2em}$e := \mathrm{fl}(b - \mathrm{fl}(s - a))$ &\\
\midrule
$(p,e):=\TwoProd(a,b)$     & 2 \\
\hspace{2em}$p := \mathrm{fl}(ab)$ & \\
\hspace{2em}$e := \mbox{FMA}(a, b, -p)$ & \\
\bottomrule
\end{tabular}}
\end{table}

%\subsection{Overall structure of the proposed FMA}
We compute $z=x\cdot y+c$ for $K$-component numbers $x,y,c$ ($K=2,3,4$ for DW/TW/QW).
The whole computation consists of three stages:
(i) product expansion (\TP{} of $x_iy_j$),
(ii) accumulation by weight (terms of the same power of two are collected by a chain
of \TS{}), and
(iii) normalization (reshaping into $K$ non-overlapping components).
Products $x_iy_j$ with $i+j\ge K$ (low-order terms of weight $O(u^K)$ and below) are
never formed in the first place, and the error bound of this truncation is
machine-proved in Section~\ref{sec:verify}.
Algorithm~\ref{alg:dwfma} shows DW, Algorithm~\ref{alg:twfma} TW, and
~\ref{alg:qwfma} QW.

\subsection{DW-FMA (proposed)}
\begin{algorithm}[h]
\caption{DW-FMA $z=x\cdot y+c$ ($K=2$, 17 flops)}\label{alg:dwfma}
\Input{$x=(x_0,x_1),\ y=(y_0,y_1),\ c=(c_0,c_1)$ (each non-overlapping)}
\Output{$z=(z_0,z_1)$}
$(P_{00},E_{00})\gets\TwoProd(x_0,y_0)$\;
$P_{01}\gets \mathrm{fl}(x_0y_1);\quad P_{10}\gets \mathrm{fl}(x_1y_0)$\tcp*{$x_1y_1$ truncated at $O(u^2)$}
$\ell\gets \mathrm{fl}(P_{01}+P_{10});\quad v\gets \mathrm{fl}(E_{00}+c_1);\quad w\gets \mathrm{fl}(v+\ell)$\;
$(s,t)\gets\TwoSum(P_{00},c_0)$\;
$t_p\gets \mathrm{fl}(t+w)$\;
$(z_0,z_1)\gets\FastTwoSum(s,t_p)$\tcp*{tail: $\mathrm{exponent}(s)\!\ge\!\mathrm{exponent}(t_p)$ is proved}
\end{algorithm}

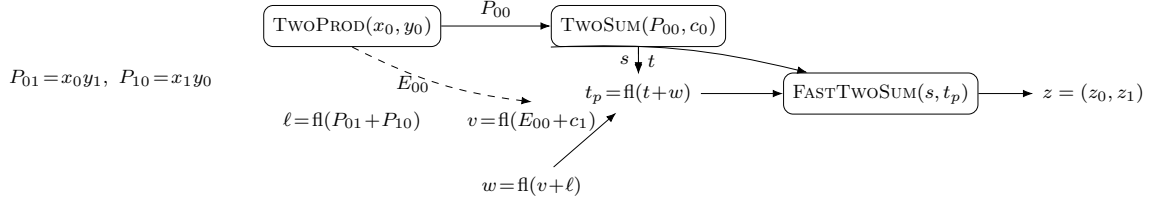
\begin{figure*}[h]\centering
\resizebox{0.96\textwidth}{!}{%
\input{fig/dwnet_en}}
\caption{Data flow of DW-FMA (solid: high-order word, dashed: error word).
The tail \FTS{} performs the non-overlapping normalization.}\label{fig:dwnet}
\end{figure*}

The data flow of Algorithm~\ref{alg:dwfma} is shown in Fig.~\ref{fig:dwnet}.
Each \TP{}/\TS{}/\FTS{} is one gate; solid lines carry the high-order word and dashed
lines the error word (the low-order word).
The truncated term $x_1y_1$ is not drawn ($O(u^2)$).
The output is made non-overlapping by the single tail \FTS{} gate, which reduces the
flop count.

\subsubsection{The div/sqrt-safe version ($20$ flops)}
The premise $\mathrm{exponent}(s)\ge\mathrm{exponent}(t_p)$ of the tail \FTS{} has been
machine-proved under the assumption that the inputs $x,y,c$ are non-overlapping.
In the iterations for division and square root described later, the intermediate values
(residuals) enter the MAC while still overlapping, so this assumption does not hold;
we therefore use a div/sqrt-safe version ($17\to20$ flops) whose tail is put back to
\TS{} (Section~\ref{sec:divsqrt}).
As \TS{} is unconditionally exact, the safe version agrees bitwise with the standard
version whenever the premise holds and is safer than the standard version only when it
does not.

\subsection{TW-FMA (proposed)}
\begin{algorithm}[h]
\caption{TW-FMA $z=x\cdot y+c$ ($K=3$, 72 flops)}\label{alg:twfma}
\Input{$x,y,c\in\mathbb{R}^3$ (non-overlapping)}
\Output{$z\in\mathbb{R}^3$}
$(P_{00},E_{00})\gets\TwoProd(x_0,y_0)$; $(P_{01},E_{01})\gets\TwoProd(x_0,y_1)$; $(P_{10},E_{10})\gets\TwoProd(x_1,y_0)$\;
$P_{02}\gets \mathrm{fl}(x_0y_2);\ P_{11}\gets \mathrm{fl}(x_1y_1);\ P_{20}\gets \mathrm{fl}(x_2y_0)$\;
$\sigma\gets\mathrm{fl}(\mathrm{fl}(P_{02}+P_{20})+P_{11})$\tcp*{symmetric: transposed pair first}
$G\gets\mathrm{fl}(\mathrm{fl}(\mathrm{fl}(E_{01}+E_{10})+\sigma)+c_2)$\;
$(A,q_1)\gets\TwoSum(P_{01},P_{10})$; $(A,q_2)\gets\TwoSum(A,E_{00})$; $(A,q_3)\gets\TwoSum(A,c_1)$\;
$G\gets \mathrm{fl}(G+\mathrm{fl}(\mathrm{fl}(q_1+q_2)+q_3))$\;
$(B,r)\gets\TwoSum(P_{00},c_0)$; $(m_1,m_2)\gets\TwoSum(r,A)$; $m_2\gets \mathrm{fl}(m_2+G)$\;
\tcp{normalization: 3 passes (\FTS{} only where its premise is provable)}
$(w_0,w_1)\gets\FastTwoSum(B,m_1)$;\ $(w_1,w_2)\gets\TwoSum(w_1,m_2)$\tcp*{pass 1}
$(w_0,w_1)\gets\TwoSum(w_0,w_1)$;\ $(w_1,w_2)\gets\FastTwoSum(w_1,w_2)$\tcp*{pass 2}
$(z_0,w_1)\gets\FastTwoSum(w_0,w_1)$;\ $(z_1,z_2)\gets\FastTwoSum(w_1,w_2)$\tcp*{pass 3}
\end{algorithm}

\begin{figure*}[htb]\centering
\resizebox{\ifdim\width>0.96\textwidth 0.96\textwidth\else\width\fi}{!}{\input{fig/twnet_en}}
\caption{Data flow of TW-FMA.  Solid lines are the representatives of each level and
dashed lines the error words that drop to the level below.
Each level is collected by a \TS{} chain and the spilled error is sent downwards,
forming a triangular structure.}\label{fig:twnet}
\end{figure*}
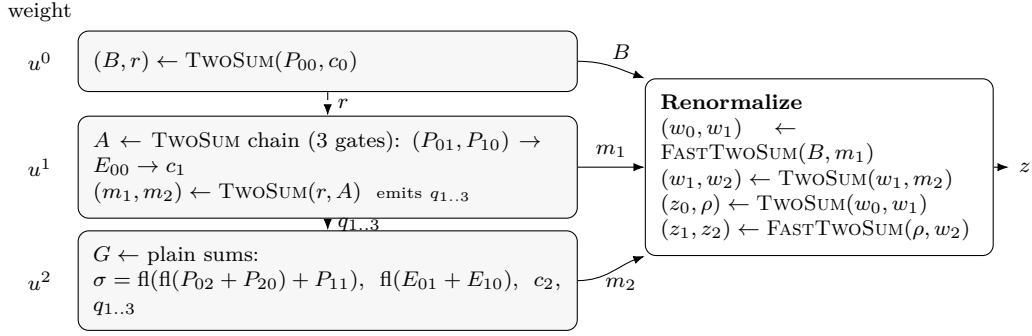

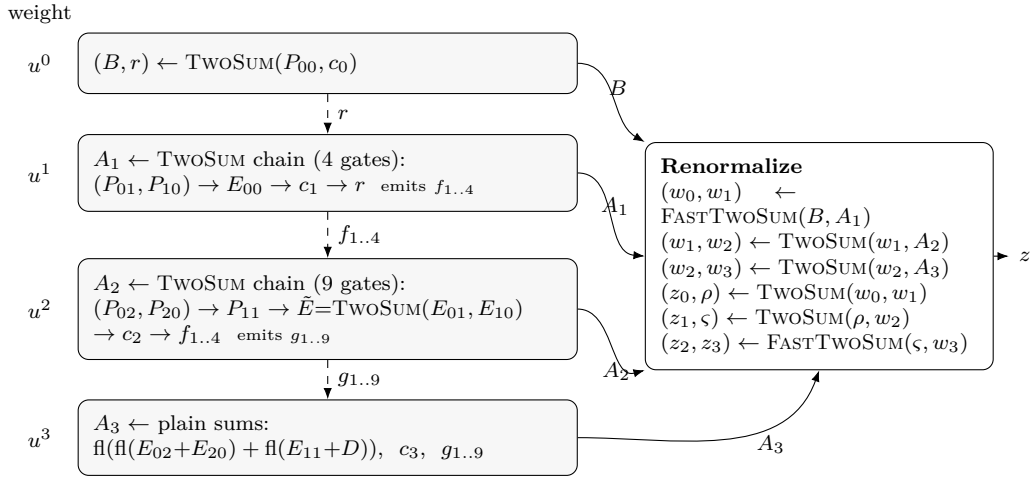
\begin{figure*}[htb]\centering
\resizebox{\ifdim\width>0.96\textwidth 0.96\textwidth\else\width\fi}{!}{\input{fig/qwnet_en}}
\caption{Data flow of QW-FMA.  Increasing $K$ by one adds one level and lengthens the
\TS{} chain of that level (3 stages in TW $\to$ 9 stages in QW).
The growth of the operation count ($72\to176$) is approximately proportional to this chain
length.}\label{fig:qwnet}
\end{figure*}

For $K\ge3$, the terms of the product expansion are naturally split into weight levels
($u^0,u^1,u^2,\dots$).
The skeleton of the proposed FMA has a triangular structure that collects each level
with a chain of \TS{} while keeping the errors, and drops the spilled errors to the
next lower level. Finally, the level representatives $B,A_1,A_2,\dots$ are made
non-overlapping by a single normalization chain, as illustrated in
Fig.~\ref{fig:twnet} (TW) and Fig.~\ref{fig:qwnet} (QW).
The essential point is that the error word of a \TS{} (dashed line) always flows into the
level one step below.  At this point the machine proof fails with two normalization
passes and becomes possible with three.

\subsubsection{The div/sqrt-safe version ($84$ flops)}
In the iterations for division and square root the inputs are not guaranteed to be
non-overlapping expansions, so the input assumption of this paper does not hold and the
proofs of Section~\ref{sec:verify} do not apply as they are.  For division and square
root we therefore use a conservative safe version ($84$ flops) in which all four \FTS{}
gates are replaced by \TS{} while the normalization keeps the same 3 passes as the
standard version (Section~\ref{sec:divsqrt}).

\subsection{QW-FMA (proposed)}
The proposed FMA algorithm QW is shown in Algorithm~\ref{alg:qwfma}; part of the
computation is abbreviated.
Here too, as in TW, the terms of the product expansion split naturally into weight
levels ($u^0,u^1,u^2,\dots$), each level is collected by a \TS{} chain that keeps the
errors, and the spilled errors drop to the level one step below, forming a triangular
structure; the machine proof fails with two normalization passes and becomes possible
with five.
\begin{algorithm}[h]
\caption{QW-FMA $z=x\cdot y+c$ ($K=4$, 176 flops)}\label{alg:qwfma}
\Input{$x,y,c\in\mathbb{R}^4$ (non-overlapping)}
\Output{$z\in\mathbb{R}^4$}
\tcp{expansion: \TP{} for $i\!+\!j\le2$, plain product $D$ for $i\!+\!j\!=\!3$, truncation for $i\!+\!j\!\ge\!4$ ($O(u^4)$)}
$(P_{ij},E_{ij})\gets\TwoProd(x_i,y_j)\ (i\!+\!j\!\le\!2)$;
$D\gets\mathrm{fl}(\mathrm{fl}(x_0y_3+x_3y_0)+\mathrm{fl}(x_1y_2+x_2y_1))$\tcp*{symmetric: transposed pair first}
$(B,r)\gets\TwoSum(P_{00},c_0)$\;
\tcp{L1: accumulation of weight $u$ ($A_1$ and the errors $f_{1..4}$)}
$(A_1,f_1)\gets\TwoSum(P_{01},P_{10})$;\ $(A_1,f_2)\gets\TwoSum(A_1,E_{00})$;\ $(A_1,f_3)\gets\TwoSum(A_1,c_1)$;\ $(A_1,f_4)\gets\TwoSum(A_1,r)$\;
\tcp{L2: 9-fold \TS{} accumulation of weight $u^2$ ($A_2$ and the errors $g_{1..9}$)}
$(A_2,g_1)\gets\TwoSum(P_{02},P_{20})$;\ $(A_2,g_2)\gets\TwoSum(A_2,P_{11})$\;
$(\tilde{E},g_4)\gets\TwoSum(E_{01},E_{10})$;\ $(A_2,g_3)\gets\TwoSum(A_2,\tilde{E})$\tcp*{symmetric: transposed pair first}
$(A_2,g_5)\gets\TwoSum(A_2,c_2)$;\ \dots\ ;\ $(A_2,g_9)\gets\TwoSum(A_2,f_4)$\;
\tcp{L3: rounded sum of weight $u^3$}
$A_3\gets\mathrm{fl}(\mathrm{fl}(\mathrm{fl}(\mathrm{fl}(E_{02}+E_{20})+\mathrm{fl}(E_{11}+D))+c_3)+(\text{tree sum of }g_{1..9}))$\tcp*{symmetric: transposed pair first}
\tcp{normalization: 5 passes (\FTS{} only where its premise is provable)}
$(w_0,w_1)\gets\FastTwoSum(B,A_1)$;\ $(w_1,w_2)\gets\TwoSum(w_1,A_2)$;\ $(w_2,w_3)\gets\TwoSum(w_2,A_3)$\tcp*{pass 1}
$(w_0,w_1)\gets\TwoSum(w_0,w_1)$;\ $(w_1,w_2)\gets\TwoSum(w_1,w_2)$;\ $(w_2,w_3)\gets\FastTwoSum(w_2,w_3)$\tcp*{pass 2}
$(w_0,w_1)\gets\TwoSum(w_0,w_1)$;\ $(w_1,w_2)\gets\FastTwoSum(w_1,w_2)$;\ $(w_2,w_3)\gets\FastTwoSum(w_2,w_3)$\tcp*{pass 3}
$(w_0,w_1)\gets\FastTwoSum(w_0,w_1)$;\ $(w_1,w_2)\gets\FastTwoSum(w_1,w_2)$;\ $(w_2,w_3)\gets\FastTwoSum(w_2,w_3)$\tcp*{pass 4}
$(z_0,w_1)\gets\FastTwoSum(w_0,w_1)$;\ $(z_1,w_2)\gets\FastTwoSum(w_1,w_2)$;\ $(z_2,z_3)\gets\FastTwoSum(w_2,w_3)$\tcp*{pass 5}
\end{algorithm}

\subsubsection{The div/sqrt-safe version ($206$ flops)}
For the same reason as in TW, in the div/sqrt context no \FTS{} is used at all
(all \TS{}) and the number of passes is kept at 5 as in the standard version;
the renormalization costs $90$ flops, $206$ flops in total.
For DW/TW/QW alike the safe version is heavier than the standard one
($20\!>\!17$, $84\!>\!72$, $206\!>\!176$).

\subsection{Design of the normalization: gate choice and number of passes}\label{sec:renorm}
The normalization stage of the proposed FMA is structured as one \textbf{cascade sweep}
of \TS{}/\FTS{} over adjacent pairs forming one pass, repeated several times.
The design freedom consists of two choices: whether each gate can be an \FTS{}, and how
many passes to repeat.

\paragraph*{(a) Gate choice.}
Each gate is chosen according to whether its \FTS{} premise can be machine-proved for
every precision $p\ge8$.
The \FTS{} validity condition of FPANVerifier is
$\mathrm{exponent}(a)\ge\mathrm{exponent}(b)$ (or one of the operands being $0$),
not the textbook $|a|\ge|b|$.
Using \FTS{} only at provable positions lowers the flop count while preserving the
property that \FTS{} is used but all of its premises are proved.
Since each replacement changes the network, we iterated to a \textbf{fixed point} at
which the verifier stops suggesting new \FTS{} candidates, re-verifying the output
guarantees at each step.
As a result, \FTS{} can be used in 1 of 1 gate for DW, 4 of 6 for TW and 10 of 15 for
QW.

\paragraph*{(b) Number of passes.}
\textbf{Gate choice alone cannot guarantee the non-overlapping property of the output.}
If the number of passes is insufficient, the non-overlapping property cannot be proved
no matter which \FTS{} is demoted to \TS{}.
Indeed, for QW with 2 normalization passes, even after replacing every \FTS{} by \TS{},
the weakest relation \texttt{p\_dominates}
($\mathrm{exponent}(z_i)\ge\mathrm{exponent}(z_{i+1})+p$) is refuted for the first two
pairs.
Each additional pass fixes one more pair from the top, and with \textbf{3 passes for TW
and 5 passes for QW} every pair is proved with the strongest
\texttt{strongly\_dominates} (Section~\ref{sec:verify}).

The normalization stages of Algorithms~\ref{alg:twfma} and \ref{alg:qwfma} are the
result of combining these two choices; they extend the originally devised algorithms
(2 passes) by appending passes at the tail (TW: $66\to72$, QW: $146\to176$ flops).

\subsection{The Exact FMA (distillation version): algorithm and error bound}\label{sec:exact}
As a reference implementation paired with the proposed FMA we use an Exact FMA
(distillation FMA) with a binomial-type error bound in which the result-relative term
dominates within a certain cancellation range.
It is a yardstick for the attainable accuracy rather than the main proposal of this
paper, but the contrast is essential for understanding the error semantics of the
proposed FMA.

Distillation \cite{ORO05} refers to the following operation on an expansion vector:
$v=(v_0,\dots,v_{M-1})$ ordered by weight: pass adjacent pairs of words through \TS{}
from the low end to the high end, such that the overflow is carried upward, thereby
bringing the whole vector into an almost non-overlapping form (each word is at most $u$
times the word one position above it).
One pass from the low end to the high end is called a sweep, and the more sweeps are
repeated, the closer the expansion becomes to non-overlapping (much as bubble sort
places one more element in order per pass).

\begin{algorithm}[h]
\caption{Exact FMA $z=x\cdot y+c$ ($K$ components, $n=K{+}1$ sweeps, 194/552/1178 flops)}\label{alg:exfma}
\Input{$x,y,c\in\mathbb{R}^K$ (non-overlapping)}
\Output{$z\in\mathbb{R}^K$}
\tcp{(1) error-free expansion: $M=2K^2+K$ terms (DW 10 / TW 21 / QW 36)}
\lForEach{$0\le i,j<K$}{$(P_{ij},E_{ij})\gets\TwoProd(x_i,y_j)$}
$v \gets \{P_{ij}\}\cup\{E_{ij}\}\cup\{c_0,\dots,c_{K-1}\}$\tcp*{$\textstyle\sum v = x\cdot y+c$ holds exactly}
\tcp{(2) symmetrization: combine the transposed pair first by \TS{} (for commutativity)}
\lForEach{$i<j$}{apply $(v_{ij},v_{ji})\gets\TwoSum(v_{ij},v_{ji})$ on both the $P$ and the $E$ side}
\tcp{(3) distillation: $n=K{+}1$ low-to-high \TS{} sweeps (DW $3$ / TW $4$ / QW $5$)}
\For{$p\gets1$ \KwTo $n$}{
  \lFor{$i\gets M-2$ \KwTo $0$}{$(v_i,v_{i+1})\gets\TwoSum(v_i,v_{i+1})$}
}
\tcp{finishing the leading $K$ words}
\lFor{$i\gets K-1$ \KwTo $0$}{$(v_i,v_{i+1})\gets\TwoSum(v_i,v_{i+1})$}
$z\gets(v_0,\dots,v_{K-1})$\tcp*{$v_K,\dots,v_{M-1}$ are truncated}
\end{algorithm}

%\paragraph{The only source of error is the truncation.}
The decisive difference from the proposed FMA is where information is discarded.
\begin{itemize}
\item \textbf{Proposed FMA}: the higher-order products ($i+j\ge K$) are never formed,
      and the diagonal sums and error sums are rounded by plain additions.
      As the discarded amount is determined by the magnitudes of $|xy|$ and $|c|$,
      the error bound is relative to $(|xy|+|c|)$
      ($35u^2/187u^3/822u^4$, Section~\ref{sec:verify}).
\item \textbf{Exact FMA}: the expansion in step (1) is exact.
      As $x\cdot y=\sum_{i,j}x_iy_j$ and each $x_iy_j$ is split exactly into
      $P_{ij}+E_{ij}$ by \TP{}, the identity $\textstyle\sum v = x\cdot y + c$ holds
      without any rounding error.
      All the \TS{} in steps (2) and (3) are EFTs as well, and therefore the sweeps preserve the
      multiset sum. Hence, the only error is the final truncation of
      $v_K,\dots,v_{M-1}$:
      \[
        z-(x\cdot y+c) \;=\; -\sum_{i\ge K} v_i .
      \]
\end{itemize}
Because the sweeps bring the expansion into an almost non-overlapping form
($|v_{i+1}|\lesssim u|v_i|$), the discarded tail is determined by $v_0$, that is by the
magnitude of the result; in the idealized limit where complete non-overlapping is
achieved, writing $\tau:=x\cdot y+c$ we obtain the result-relative form
$|z-\tau|\le C_K u^K|\tau|$.
Here $C_K$ is a constant that originates from the fact that the number of sweeps $n$ is
finite (with an infinite number of sweeps, the expansion converges to Priest's
non-overlapping expansion\cite{Priest91} and $C_K\to O(1)$).

%\paragraph{The rule $n=K{+}1$ for the number of sweeps (a correction made in this paper).}\label{par:sweeprule}
The consequence of the finiteness of the sweeps is not merely that the constant $C_K$
becomes large. The tail that remains after $n$ sweeps contains, in addition to the
result-relative term, an input-relative residual, such that the bound takes the two-term
form
\[
  |z-\tau| \;\le\; P_K(u)\,u^K\,|\tau| \;+\; Q_K(u)\,u^{\,n}\,\bigl(|x\cdot y|+|c|\bigr)
\]
(Appendix~\ref{sec:proofex}). If $n=K$ is taken, the second term becomes
$u^K(|xy|+|c|)$, of the same order as the first, and under deep cancellation
($|\tau|\ll|xy|+|c|$), the second term dominates and the result-relative property is lost.
DW had $n=3=K{+}1$ from the beginning for its $M=10$ terms, whereas TW and QW used
$n=K$ ($3$ and $4$) and, in this sense, did not justify the name ``Exact''.

In this study, we unify the number of sweeps to $n=K{+}1$ (TW $3\to4$, QW $4\to5$; DW
unchanged). This is not a new design but a unification of the rule, and the only change
in the implementation is the trip count of the distillation loop.
As the symmetrization stage is untouched, commutativity is preserved (the expansion
vector after symmetrization is position-wise invariant under the transposition
$x\leftrightarrow y$, and the additional sweep is a deterministic map with a fixed order).

The condition under which the result-relative reading dominates,
$|\tau|\gtrsim(Q_K/P_K)\,u\,(|xy|+|c|)$, permits cancellation down to a depth of about
$40.4$ bits for TW and about $37.9$ bits for QW at $p=53$.
This is on a par with DW (about $43.9$ bits), so within this cancellation range the
three word lengths DW/TW/QW all share the same semantics of being result-relative at
the order $u^K$.

%\paragraph{What being result-relative means.}
Under the strong cancellation $c\approx-xy$, we have $|\tau|\ll|xy|+|c|$.
In this case, the relative error of the proposed FMA, whose bound is relative to
$(|xy|+|c|)$, can deteriorate with respect to the result, whereas the fully distilled
Exact FMA is unaffected down to a cancellation of approximately $38$--$44$ bits.
The measurements also show that the Exact FMA maintains its advantage under cancellation
(the division of Section~\ref{sec:divsqrt}).
Its operation count, however, is $194/552/1178$, i.e.\ $6.7$--$11.4$ times that of the
proposed FMA ($17/72/176$), so it is useful only for local, accuracy-critical uses.
This is why we position the Exact FMA as a reference for the attainable accuracy and
the proposed FMA as the practical replacement.

\subsection{Theoretical operation counts}\label{sec:flops}

We summarize the theoretical operation counts of the proposed FMA algorithms.
Table~\ref{tab:flop} lists the operation count per accumulating MAC
$z=a\cdot b+c$.

\begin{table*}[htb]\centering
\caption{Theoretical operation count per accumulating MAC $z=a\cdot b+c$ (a hardware FMA counts as $1$)}\label{tab:flop}
\setlength{\tabcolsep}{2.5pt}
\footnotesize
\begin{tabular}{@{}lrrrr l@{}}
\toprule
variant & DW (2) & TW (3) & QW (4) & BF/FMA & remarks\\
\midrule
Q (existing mul+add) & $20^{\dagger}$ & $\sim$125\,+renorm$^{\ddagger}$ & $\sim$162\,+renorm4$^{\ddagger}$ & --- & \\
BF (mul\_bf+add\_bf) & 29 & 96 & 209 & --- & branch-free\\
\textbf{FMA (proposed, proved)} & \textbf{17}$^{\S}$ & \textbf{72}$^{\S}$ & \textbf{176}$^{\S}$ & --- & non-overlapping proved\\
FMA (2-pass normalization, reference) & 17 & 66 & 146 & --- & non-overlapping not provable$^{\P}$\\
FMA (div/sqrt-safe) & 20 & 84 & 206 & --- & all \TS{}; for different input assumptions\\
Exact FMA ($n{=}K{+}1$) & 194 & 552 & 1178 & --- & binomial-type bound, commutative\\
\midrule
flop ratio BF/FMA (proved) & 1.71 & 1.33 & 1.19 & & \\
\bottomrule
\end{tabular}

\smallskip
{\footnotesize
$\dagger$ DW-Q (sloppy) is branch-free and has a fixed count.\quad
$\ddagger$ The normalizations renorm/renorm4 of TW/QW-Q contain value-dependent
branches (executed as scalar code per lane) and do not have a fixed operation count.
This is the main reason why Q is slow, and it is precisely what BF and the proposed FMA
eliminate.\\
$\S$ DW 17 / TW 72 / QW 176 is the configuration for which the anchor-relative error
bounds, the FastTwoSum premises and \textbf{the non-overlapping property of the output}
have been machine-proved by FPANVerifier (Section~\ref{sec:verify}).
All timing measurements in this paper use this proved version.\\
$\P$ TW 66 / QW 146 is the earlier version with 2 normalization passes; its error
bounds and \FTS{} premises are likewise proved, but the non-overlapping property of the
output cannot be proved (Section~\ref{sec:nonoverlap}).  Measured values are also given
for comparison.}
\end{table*}

The flop ratio BF/FMA of the proved version is $1.71/1.33/1.19$.
The ratio approaches 1 as $K$ grows because the number of normalization passes required
to guarantee the non-overlapping property increases with $K$ as $1/3/5$, so the
normalization takes a growing share of the whole computation.
If a kernel is sufficiently compute-bound, this theoretical ratio appears directly in
the measurements: in the AVX2 GEMM (DD \spd{1.51}, Table~\ref{tab:time_x86}) and in the
CUDA GEMM/GEMV on GB10 (DD \spd{1.71}, Table~\ref{tab:cuda}) the type with the largest
operation-count ratio, DW, shows the largest speed-up.
The Arm/NEON GEMM, on the other hand, shows little dependence on the type,
\spd{1.40}--\spd{1.46} (Table~\ref{tab:time_arm}); as discussed later this is decided by
whether the kernel is compute-bound or bandwidth-bound.
The Exact FMA needs $11.4/7.7/6.7$ times as many operations as the proposed FMA, and it
is measured to be about 4--12 times slower in the three NEON kernels.

\subsection{Machine proof of the error bounds}\label{sec:verify}
The error bounds of the proposed FMA and the validity of its normalization gates were
machine-proved with FPANVerifier\cite{ZhangAiken25CAV} of Zhang and Aiken, an abstract
interpreter for floating-point arithmetic that uses the SMT solver Z3.
Since the verification description \texttt{.fpan} treats the precision $p$ symbolically
(assuming $p\ge8$), the proofs hold at once for all formats with $p\ge8$, including
binary32/64/128.

%\subsubsection{Certified error bounds}
The best provable bound of each discarded wire is found by bisection, and their sum is
the overall bound.  The results are shown in Table~\ref{tab:certbound}.
The truncated products $x_iy_j$ ($i+j\ge K$) are certified by \texttt{bound} through
the exact identity $x_iy_j=P_{ij}+E_{ij}$ of \TP{}, covering both the high word
$P_{ij}$ and the low-order error $E_{ij}$; the higher-order products that the
implementation never forms (for QW, $x_2y_3$, $x_3y_2$ and $x_3y_3$) are declared and
accounted for in the same way.
The \texttt{bound} of the machine proof is a relative bound whose \textbf{anchor} is
the leading product $P_{00}=\mathrm{fl}(x_0y_0)$ (product side) and
$\max(|P_{00}|,|c_0|)$ (mixed side); the anchor-relative coefficient sums are the
low-order sums $34/184/812$ plus a high-order tail coming from the truncated products:
DW $(34+4u)u^2$, TW $(184+12u+4u^2)u^3$ and QW $(812+20u+12u^2+4u^3)u^4$.
To move these to the $(|xy|+|c|)$-relative form, put $S_K=u+\cdots+u^{K-1}$; the
non-overlapping condition of the inputs $|x_{i+1}|\le u|x_i|$ etc.\
(Appendix~\ref{sec:proof}) yields
$|P_{00}|\le(1+u)|x_0y_0|\le(1+u)|x\cdot y|/(1-S_K)^2$ and $|c_0|\le|c|/(1-S_K)$,
so it suffices to multiply the coefficients by $(1+u)/(1-S_K)^2$.
Rounding up at the worst case over $p\ge8$ (namely $p=8$) gives the input-relative
machine-proved bounds $35u^2$/$187u^3$/$822u^4$, valid for every $p\ge8$ including
binary32/64/128 (Table~\ref{tab:certbound}).
The verification time of FPANVerifier (z3), with a 20-core parallel sharded run on
GB10 (Table~\ref{tab:env}), is approximately 34 s, 22 min, and 1.6 h
to finish all queries of DW, TW and QW, respectively.

We also show the values obtained independently by a chain of paper-and-pencil
inequalities in the style of Joldeş--Muller--Popescu (Appendix~\ref{sec:proof}).
At practical precision $p\ge8$ the manual proof is $2.3$--$2.7$ times sharper, and the
two serve as crosschecks for each other. Conversely, at low precision, $p=2$,
FPANVerifier is the sharper of the two.

\begin{table}[h]\centering
\caption{Error bound of the proposed FMA, $|z-(xy+c)|\le C_K\,u^K\,(|xy|+|c|)$ ($u=2^{-p}$).
Neither the machine-proved nor the manual bound depends on the number of normalization
passes.}\label{tab:certbound}
\footnotesize
\setlength{\tabcolsep}{3pt}
\begin{tabular}{lrrlr}
\toprule
     & FPAN-    & manual & manual,    & measured\\
type & Verifier & proof  & asymptotic & max\\
\midrule
DW & $35\,u^2$  & $\mathbf{14}\,u^2$  & $(13|xy|+4|c|)u^2$   & $4.9\,u^2$\\
TW & $187\,u^3$ & $\mathbf{67}\,u^3$  & $(65|xy|+8|c|)u^3$   & $5.2\,u^3$\\
QW & $822\,u^4$ & $\mathbf{360}\,u^4$ & $(348|xy|+44|c|)u^4$ & $4.8\,u^4$\\
\bottomrule
\end{tabular}

\smallskip
{\footnotesize The values in the FPANVerifier column are the anchor-relative
machine-proved bounds moved to the input-relative ($(|xy|+|c|)$-relative) form by the
analytic transformation (Section~\ref{sec:verify}).
Since the FPANVerifier proof treats the precision $p$ symbolically ($p\ge8$), it holds
at once for all formats with $p\ge8$ including binary32/64/128.
``Measured max'' is the maximum of $\eta=|z-(xy+c)|/(|xy|+|c|)$, reported as
$\max\eta/u^K$, observed in the falsification test with MPFR
(400{,}000 trials per type and per backend, no violation).}
\end{table}

\subsubsection{Non-overlapping property of the output and the number of normalization passes}\label{sec:nonoverlap}
To hand the output downstream as a ``non-overlapping expansion'', the dominance
relation $z_i\vartriangleright z_{i+1}$ between adjacent components must be guaranteed.
FPANVerifier can test this with predicates of different strength (from the strongest:
\texttt{strongly\_dominates}, \texttt{qd\_dominates}, \texttt{ulp\_dominates},
\texttt{s\_dominates}, \texttt{p\_dominates}).
Table~\ref{tab:passes} shows the results of the test while varying the number of
normalization passes.

\begin{table}[h]\centering
\caption{Number of normalization passes and the strongest provable dominance relation
($-$: even the weakest \texttt{p\_dominates} is refuted)}\label{tab:passes}
\setlength{\tabcolsep}{4pt}
\small
\begin{tabular}{llll}
\toprule
type & passes & $z_0\vartriangleright z_1$ & $z_1\vartriangleright z_2$ ($z_2\vartriangleright z_3$)\\
\midrule
TW & 2 (former) & \texttt{ulp} & \texttt{strongly}\\
TW & \textbf{3 (this paper)} & \textbf{\texttt{strongly}} & \textbf{\texttt{strongly}}\\
\midrule
QW & 2 (former) & $-$ & $-$\ \ (\texttt{strongly})\\
QW & 3 & $-$ & \texttt{p}\ \ (\texttt{p})\\
QW & 4 & \texttt{strongly} & \texttt{s}\ \ (\texttt{strongly})\\
QW & \textbf{5 (this paper)} & \textbf{\texttt{strongly}} & \textbf{\texttt{strongly}}\ \ (\textbf{\texttt{strongly}})\\
\bottomrule
\end{tabular}
\end{table}

For the TW and QW FMAs, the non-overlapping property of the renormalization part cannot
be proved with 2 passes.
For TW the two strong relations of the leading pair are refuted, and for QW even the
weakest \texttt{p\_dominates} is refuted for the first two pairs (the same holds after
replacing every \FTS{} by \TS{}).
Each additional pass fixes one more pair from the top, and with 3 passes for TW and 5
passes for QW every pair is proved with \texttt{strongly\_dominates}.
The abstract counterexample at refutation is the \textbf{complete cancellation} ($B=0$)
of the $\TwoSum(P_{00},c_0)$ of $L_0$ at $x\cdot y\approx-c$; once the leading word
vanishes it has to be rebuilt from the lower words, which a fixed 2-pass scheme cannot
recover.

%\paragraph{An important caveat.}
Note that the abstraction of FPANVerifier (SELTZO) is \textbf{sound but not complete},
so REFUTED does not mean that the original 2-pass version is actually wrong.
In more than $2\times10^6$ measured trials sweeping the cancellation strength from
$2^{-20}$ to $2^{-110}$, \textbf{not a single \texttt{p\_dominates} violation was
observed even for the 2-pass version}.
What the additional passes buy is therefore a \textbf{machine-checked guarantee}, not
the fix of a known bug.
Making this distinction explicit, this paper gives priority to the guarantee and adopts
the proved version (TW $72$ / QW $176$ flops) as the default.
The error bounds of Table~\ref{tab:certbound} do not depend on the number of passes
(Appendix~\ref{sec:proof}).

%----------------------------------------------------------%
% Related work
%----------------------------------------------------------%
\section{Related work}\label{sec:related}

There have been several attempts to exploit FMA in multiple-precision arithmetic,
but the term ``FMA'' refers to different computational objects in different studies.
In this section, we consider three representative lines of work---the error of a scalar FMA by
Boldo and Muller, the DW-FMA for matrix multiplication by Ozaki and Koizumi, and the FMA
emulation by Muller et al.---and organize the differences from the FMA proposed here
(Table~\ref{tab:related}).

\begin{table*}[htb]\centering
\caption{Positioning of the related work on FMA ($u=2^{-p}$; flops are per MAC, a hardware FMA counting as $1$)}\label{tab:related}
\footnotesize
\setlength{\tabcolsep}{2pt}
\begin{tabular}{@{}p{26mm}p{40mm}p{18mm}p{44mm}p{22mm}@{}}
\toprule
study & target of the computation & flops & applicability condition & proof system\\
\midrule
\textbf{This work} & \textbf{$K$-component FMA} $x\cdot y+c$ (DW/TW/QW, $K$ components both in and out) & $17/72/176$ & non-overlapping inputs, $p\ge8$\newline (no additional dominance condition of type $|c|\ge2|ab|$; every \FTS{} premise proved) & SMT\newline (FPANVerifier)\\
\addlinespace
Zhang--Aiken \cite{ZhangAiken25SC,ZhangAiken25CAV} & branch-free \texttt{mul}$+$\texttt{add} (the BF baseline of this study) & $29/96/209$ & none & SMT\\
\addlinespace
Boldo--Muller \cite{BoldoMuller11} \textsc{ErrFma} & the \textbf{exact error} of a \textbf{scalar} FMA\newline $ax{+}y=r_1{+}r_2{+}r_3$ & $20$ & even radix, $p\ge3$ (subnormals allowed) & Coq\\
\addlinespace
Boldo--Muller \cite{BoldoMuller11} \textsc{ErrFmaAppr} & approximation of the above (error $O(u^2)|z|$) & $12$ & even radix, $p\ge3$ (subnormals allowed) & Coq\\
\addlinespace
Ozaki--Koizumi \cite{OzakiKoizumi25} \textsc{FastFma\_DW} & \textbf{DW-FMA} ($K{=}2$ only, embedded in a matrix product) & $6$ & \textbf{$c-x$ must be exact} ($|c|\ge2|ab|$ etc.; guaranteed by a shift matrix $E$ in the matrix product) & on paper\\
\addlinespace
Jeannerod et al.\ \cite{JJLM26} & \textbf{tight error bounds} for a family of fast DW-FMA kernels (not a new algorithm) & --- & \textbf{dominant regime} $|c|\ge2|ab|$ plus an overlap parameter $|x_\ell|\le k\,\mathrm{ulp}(x_h)$ & on paper\newline (with worst cases)\\
\bottomrule
\end{tabular}
\end{table*}

%----------------------------------------------------------%
% Related work 1
%----------------------------------------------------------%
\subsection{Boldo--Muller: the exact error of a scalar FMA (\textsc{ErrFma})}
Boldo and Muller showed that the error produced by one hardware FMA can always be
represented exactly by two floating-point numbers, and gave \textsc{ErrFma}
($20$ operations), which computes $r_1=\circ(ax+y)$ and $r_2,r_3$ such that
$ax+y=r_1+r_2+r_3$ holds exactly\cite{BoldoMuller11}.
They also provided a formal proof in Coq and extended the result to a general even radix,
not only radix 2, and to cases with subnormals.
An approximate version \textsc{ErrFmaAppr} that requires only $12$ operations at the price
of a slight loss in accuracy is also presented.

The difference from this work lies in the level of the object being computed.
What \textsc{ErrFma} handles is the error of one FMA on scalars
$a,x,y\in\mathbb{F}$, not an FMA between multiple-precision numbers.
The FMA proposed here, in contrast, has $K$ components (DW/TW/QW) in both its input and
its output.  Among its internal EFT primitives, \texttt{two\_prod} calls the hardware
FMA, whereas \texttt{two\_sum} uses only additions and subtractions.

There is one particularly important point of contact.
The ``fma-folded diagonal terms'' listed as future work in Section~\ref{sec:concl}
(folding $a\cdot b+c$ directly with a single fused rounding would reduce the flop
counts further) were not adopted in this paper, because a
single fused rounding cannot be expressed as a composition of
\texttt{two\_prod}$+$\texttt{two\_sum} and therefore cannot be machine-proved with
FPANVerifier.
The theorem of Boldo and Muller supplies this missing gate exactly: an exact
representation of the residual of a fused $ax+y$. If an ``FMA gate'' based on their
Coq theorem could be added to FPANVerifier, the folded version would also come
within the reach of machine proof.

A remark on the operation counts is also in order.
The $20$ operations of \textsc{ErrFma} are the price of the error of one scalar FMA, and
the $17$ flops of our DW-FMA (an entire MAC for $K{=}2$) measure a different object.
Boldo and Muller themselves write that they ``could not find an algorithm with fewer
operations''; just as our $17/72/176$ is minimal within a given normalization structure,
both claims are only claims of ``the smallest known so far''.

%----------------------------------------------------------%
% Related work 2
%----------------------------------------------------------%
\subsection{Ozaki--Koizumi: a DW-FMA for matrix multiplication (\textsc{FastFma\_DW})}
Ozaki and Koizumi conducted the first error analysis of \textsc{FastTwoFma}
($x=\texttt{fma}(a,b,c)$ and $y=\texttt{fma}(a,b,c-x)$; Koizumi et al.\cite{Koizumi23}),
an approximate EFT for $ab+c$ built from only two FMAs and proposed a DW-FMA
\textsc{FastFma\_DW} based on it (only $6$ operations) along with a matrix
multiplication \textsc{Sesqui\_Word\_MM} that is embedded within it\cite{OzakiKoizumi25}.
On Ryzen 9 7950X, it remained within \spd{4.1}--\spd{4.7} of OpenBLAS \texttt{dgemm},
which is faster than pair arithmetic (PA, \spd{6.5}--\spd{7.4}) and DW
(\spd{11}--\spd{13}).

The applicability condition is decisive.
The error in \textsc{FastTwoFma} remains within $O(u^2)$ only when $c-x$ is exactly
representable; otherwise, a term $O(u)|ab|$ remains.
Specifically, when $\mathrm{sign}(c)\ne\mathrm{sign}(ab)$ and $|c|\approx|ab|$ (strong
cancellation), the relative error can deteriorate significantly.
They prove sufficient conditions, such as $|c|\ge2|ab|$, and in the matrix product, they
avoid this problem by adding a shift matrix $E$ (whose $e_{ij}$ are formed from the row
and column norms rounded upward) such that the condition holds for every $k$.

There are three differences from this work.

\textbf{(1) Presence of a premise.}  Their DW-FMA is a fused MAC that imposes conditions
on the operands, and is complete as a matrix-multiplication algorithm that comes with
the machinery to guarantee those conditions (precomputing the shift matrix $E$).
The FMA proposed here, in contrast, requires no additional dominance condition between
the operands such as $|c|\ge2|ab|$, under the common assumption that the inputs are
non-overlapping expansions; every \FTS{} premise has been machine-proved for all
precisions $p\ge8$ (Section~\ref{sec:verify}).
It can therefore be used as a drop-in replacement for \texttt{mul}$+$\texttt{add} in
AXPY/GEMV/GEMM, and it applies to division and square root through a context-dependent
safe version in which the \FTS{} are replaced by \TS{}.
This difference is not academic: as seen in Section~\ref{sec:divsqrt}, the inner MAC of
division and square root receives residuals whose non-overlapping property is not
guaranteed, so even in this work every \FTS{} of the gate placement proved for the
standard MAC had to be replaced by \TS{} (the div/sqrt-safe versions,
DW $20$/TW $84$/QW $206$).
The lesson of this paper---that a proved gate placement must be re-verified for each
context of use---explains precisely why Ozaki and Koizumi engineered the condition with
$E$ only within the context of a matrix product.

\textbf{(2) Accuracy target.}  They deliberately choose an accuracy below DW
(about $1.5$ words, ``sesqui-word'') and buy speed with it.
Indeed, at $n{=}1000$ the maximum relative error of their method is \ee{5.5}{-22}
against PA \ee{2.4}{-24} and DW \ee{9.2}{-25}, i.e.\ 2--3 digits worse\cite{OzakiKoizumi25}.
The FMA proposed here requires the accuracy of DW/TW/QW to be preserved; we
established the error bounds $35u^2/187u^3/822u^4$ by the machine proof of the
anchor-relative bounds and the analytic transformation, and measured that the rounding
agrees bitwise across five backends.

\textbf{(3) Number of components.}  They treat only $K{=}2$ (DW).  This work treats
$K{=}2,3,4$ (DW/TW/QW) uniformly, and the normalization stage can be reduced for every
$K$ (against BF, DD \spd{1.35}--\spd{1.51} and QD \spd{1.07}--\spd{1.21} in GEMM).

%----------------------------------------------------------%
% Related work 3
%----------------------------------------------------------%
\subsection{Jeannerod--Joldeş--Louvet--Muller: tight bounds with a parameterized overlap}

The condition required by \textsc{FastFma\_DW} of Ozaki and Koizumi---$|c|\ge2|ab|$,
which guarantees the exactness of $c-x$---is not an accidental trick but a structure
common to fast FMA kernels of this type.
Jeannerod, Joldeş, Louvet, and Muller have recently addressed this issue and analyzed
error bounds of a family of fast FMA kernels that approximate $ab+c$ in a DW
representation under the condition that parameterizes the overlap between high-order
and lower-order words, $|x_\ell|\le k\,\mathrm{ulp}(x_h)$ (including two- and
three-parameter versions for several DW inputs)\cite{JJLM26}.
Targeting the same dominant regime $|c|\ge2|ab|$ as Ozaki--Koizumi and its DW versions
($|c_h|\ge2|ab|$, $|c_h|\ge2|ab_h|$, and so on, which arise in the polynomial evaluation
of elementary functions and in constructions without cancellation), they give explicit
worst-case constants for each FP/DW and DW/DW variant and provide examples that
attain them, thereby tightening existing analyses (ARITH 2026).

There are two relations to this work.

\textbf{(1) Confirmation that the guarantee is conditional.}
Their results quantitatively confirm the point made in this section, namely that the
fast FMAs of the \textsc{FastFma\_DW} family have guarantees that depend on the dominant
regime and on an overlap condition.
The significance of the proposed FMA requiring no dominance condition beyond the common
assumption of non-overlapping inputs (Section~\ref{sec:verify}) becomes clearer in this
contrast.

\textbf{(2) Implications for the premises of this work itself.}
Our proof, on the other hand, is not free of an overlap assumption either: the
\texttt{.fpan} description \texttt{assume}s \texttt{strongly\_dominates}
(non-overlapping) for the input components, which corresponds to taking their parameter
$k$ in its strongest form.
In fact, for the QW 2-pass version of our original design the non-overlapping property
of the output could not be proved on FPANVerifier for inputs involving extreme
cancellation ($c\approx-xy$) (Section~\ref{sec:verify}, ``non-overlapping property of
the QW output''), and there was a risk that feeding such an output back as an input in
the recursive accumulation of an inner product would violate the \texttt{assume}.
In the 5-pass version of this paper \texttt{strongly\_dominates} is proved for every
adjacent pair of outputs, so this concern is resolved.
If, however, the non-overlapping property of the inputs themselves breaks down, i.e.\
to give a proved bound for a weakened $k$, their parameterization becomes necessary.
This is the most natural next step directly following this work, and corresponds to
relaxing the \texttt{assume} of FPANVerifier by the parameter $k$ and proving it again.

To summarize, the two occupy different points on the accuracy--speed trade-off curve.
If an accuracy of about $1.5$ words suffices and the computation can be restricted to a
matrix product, the $6$-operation \textsc{FastFma\_DW} is superior; if a general-purpose
FMA that preserves the full DW/TW/QW accuracy without additional dominance conditions
under the common assumption of non-overlapping inputs is needed, our $17/72/176$ is the
answer.
The two are not mutually exclusive: a combination such as ``\textsc{FastFma\_DW} for a
GEMM to which the $E$ shift can be applied, and the proposed FMA for division, square
root, and any cases where the condition cannot be established'' is quite natural.

%----------------------------------------------------------%
% Related work 4
%----------------------------------------------------------%
\subsection{Muller et al.: FMA emulation and the building blocks of DW/TW}
Graillat and Muller proposed algorithms that emulate the FMA and the correctly rounded sum
of three numbers (ADD3) using only round-to-nearest addition, multiplication, and
comparison, for environments without a hardware FMA\cite{GraillatMuller25}
(earlier work by Boldo and Melquiond employed round-to-odd\cite{BoldoMelquiond08}).
This is exactly the opposite of our problem setting.
This study assumes that all five evaluation targets (AVX2 / AVX-512 / NEON / SVE2 / CUDA)
have a hardware FMA and actively exploits the fact that
TwoProd$(a,b)=(\mathrm{fl}(ab),\ \mathrm{FMA}(a,b,-\mathrm{fl}(ab)))$ can be obtained in
2 operations.
In an environment without an FMA, TwoProd degenerates into Dekker's
splitting\cite{Dekker71} ($17$ operations), so the flop table of this paper
(Table~\ref{tab:flop}) itself ceases to be a valid cost model; this deserves attention.
Conversely, the emulation of Graillat and Muller becomes necessary in such an
environment, but its cost easily exceeds the savings of the proposed FMA
($29\to17$ flops and so on).

The line of work of Muller et al.\ also includes Joldeş--Muller--Popescu\cite{JMP17},
who provided tight error bounds for the building blocks of DW arithmetic.
Muller--Rideau\cite{MullerRideau22} formalized these in Coq and
Fabiano--Muller--Picot\cite{FMP19} proposed algorithms for TW (triple-word) arithmetic.
These studies established the grounds of DW/TW/QW on which this study is based; however, they
all target \texttt{mul} and \texttt{add} and do not define an FMA as a single fused
operation.
The contribution of this study is the placement of a branch-free, machine-proved
multiple-precision FMA in that gap, folding the two existing stages
\texttt{mul}$+$\texttt{add} into one and thereby achieving an actual speed-up.

%--------------------------------------%
% Benchmark tests
%--------------------------------------%
\section{Benchmark tests}\label{sec:bench}

In this section, we evaluate the performance of the proposed FMA for all six types
(binary32-based DS/TS/QS and binary64-based DD/TD/QD) in the two computing environments
(Table~\ref{tab:env}).
Section~\ref{sec:divsqrt} evaluates the division and square root and
Section~\ref{sec:blas} presents the basic linear algebra kernels AXPY/GEMV/GEMM
(serial, OpenMP and GPU parallel).

\begin{table*}[htb]\centering
\caption{Computing environments used for the benchmark tests}\label{tab:env}
\footnotesize
\setlength{\tabcolsep}{3pt}
\begin{tabular}{@{}lp{58mm}p{56mm}@{}}
\toprule
 & \textbf{GB10} (Arm) & \textbf{H100 node} (x86)\\
\midrule
CPU & Arm Cortex-X925 3.9\,GHz\newline (big.LITTLE, 10$+$10 cores) & Intel Xeon Gold 6526Y\newline (32 threads)\\
SIMD & NEON (\texttt{float64x2\_t} / \texttt{float32x4\_t})\newline SVE2 (VL$=128$ bit) & AVX2 (\texttt{\_\_m256d})\newline AVX-512 (\texttt{\_\_m512d})\\
GPU & NVIDIA GB10 (Grace--Blackwell)\newline sm\_121, 48 SM, LPDDR5X approx.\ 273 GB/s & NVIDIA H100 NVL (Hopper)\newline sm\_90, 132 SM, HBM3 approx.\ 3900 GB/s\\
\bottomrule
\end{tabular}

\smallskip
{\footnotesize The compilers are gcc 13.3 and nvcc (CUDA 13.0), with the optimization
flags \texttt{-O3 -funroll-loops}.
nvcc requires \texttt{--fmad=false -ftz=false} (indispensable for the validity of the EFTs).
The SIMD kernels of AXPY/GEMV/GEMM write the EFTs directly with intrinsics and leave no
room for contraction, but for division and square root (\texttt{divsqrt.hpp}), which
contain plain C++ expressions, \texttt{-ffp-contract=off} (and
\texttt{\#pragma STDC FP\_CONTRACT OFF}) is specified on the host side as well.
The accuracy reference is MPFR with 600 bits (DS$\approx$48 / TS$\approx$72 /
QS$\approx$96 / DD$\approx$106 / TD$\approx$159 / QD$\approx$212 bits).
Serial measurements are pinned to one performance core with \texttt{taskset} and take
the best of the three.}
\end{table*}

\begin{table*}[htb]\centering
\caption{Four variants to be compared}\label{tab:variants}
\setlength{\tabcolsep}{2pt}
\footnotesize
\begin{tabular}{ll}
\toprule
symbol & description \\
\midrule
Q        & the default \texttt{mul}+\texttt{add} (DD=sloppy, TD/QD=vec\_sum/vseb of Bailey\cite{HidaLiBailey}, with branches)\\
BF       & the built-in branch-free \texttt{mul\_bf}+\texttt{add\_bf} of Zhang--Aiken\cite{ZhangAiken25SC} (the fair baseline)\\
FMA      & the proposed branch-free FMA (\texttt{fma\_ref.c} translated to each ISA, $17/72/176$ flops)\\
Exact FMA& the result-relative distillation FMA (\texttt{fma\_exact\_ref.c}, fully distilled, $n{=}K{+}1$, $194/552/1178$ flops)\\
\bottomrule
\end{tabular}
\end{table*}

%\paragraph{Implementation.}
The scalar TwoSum / TwoProd / FastTwoSum are mapped one-to-one onto the inline
functions of BNCmatmul that use the SIMD intrinsics of each ISA
(\texttt{\_bncavx2\_*} / \texttt{\_bncavx512\_*} / \texttt{\_bncneon\_*} / \texttt{\_bncsve2\_*},
all of which use the hardware FMA), and the plain $+,\times$ are the addition and
multiplication instructions of each ISA.
The data are stored as a column-major SoA holding a \texttt{double*} (or \texttt{float*})
per component, and one instruction processes 2--8 output elements in parallel
(AVX2 4 / AVX-512 8 / NEON 2 (4 in single precision) / SVE2 2).
CUDA is scalar per thread, and the component type $R\in\{\texttt{float},\texttt{double}\}$
and the number of components $K\in\{2,3,4\}$ are implemented uniformly with C++ templates.
The four accumulating operations $z \leftarrow x\cdot y + c$ to be compared are listed in
Table~\ref{tab:variants}.

%\paragraph{On the version of the Exact FMA (important).}
The Exact FMA used in this study is the fully distilled version of Section~\ref{sec:exact}
($n=K{+}1$ sweeps).
All the numbers of the GB10 and H100 environments
(AVX2 / AVX-512 / H100 GPU) were measured using this algorithm.

\subsection{Measurement method}\label{sec:method}
Here, we describe how the numbers used in this study were produced.
All file names below are those inside the benchmark suite (\texttt{bench\_fma}).

%\paragraph{Structure of the code.}
The prototypes of the proposed FMA and the Exact FMA were written in C, and the C code
optimized with the SIMD intrinsics and CUDA functions of each target is used as the
comparison routine of the benchmark; this is the body of the operation that is actually
measured. The three kernels AXPY/GEMV/GEMM use C/C++ code that embeds them.
The division and square root use a backend-independent CUDA kernel.

%\paragraph{How accuracy is measured.}
All accuracy references were computed using MPFR at 600 bits.
For AXPY/GEMV/GEMM, we take $|\hat{r}_i-r_i|/|r_i|$ over all the elements of the output
vector or matrix, and report its maximum and mean (Table~\ref{tab:acc}).
The sizes for the accuracy measurements were determined based on the cost of the MPFR reference:
AXPY $n{=}4000$, GEMV $n{=}300$, and GEMM $n{=}128$ (note that they differ from
the sizes used for timing). The division and square root use $n{=}8192$ elements, and their
relative errors are listed in Table~\ref{tab:dsqacc}.
The ill-conditioned inner products are produced with a GEMV of $n{=}512$ at the target
condition numbers $2\times10^{3},2\times10^{8},2\times10^{13}$.
The actual numbers assigned to each kernel are described in Section~\ref{sec:testdata}.

%\paragraph{How speed is measured (CPU).}
The wall-clock time is measured with \texttt{clock\_gettime(CLOCK\_MONOTONIC)}.
After one warm-up run, the repetition count is doubled automatically until the total
measured time reaches at least $0.3$ seconds, and this block of repetitions is executed
three times, of which the minimum is taken (best of 3).
What appears in the tables is this value divided by the repetition count, i.e.\ the
number of seconds per call.
Throughout the paper, a reported FMA/BF speed-up denotes the performance ratio,
equivalently $t_{\mathrm{BF}}/t_{\mathrm{FMA}}$ when measured by execution time; hence a
value greater than 1 means that the proposed FMA is faster.
The sizes for timing are AXPY $n{=}10^6$, GEMV $n{=}2048$, GEMM $n{=}512$ and
$n{=}131072$ for division and square root; the serial measurements are pinned to one
performance core using \texttt{taskset}.
Only OpenMP measurements (Section~\ref{sec:omp}) were built using
\texttt{-DFAST\_TIME}, skipping the best-of-3 and taking a single run after the
automatic adjustment. Thread placement was fixed with
\texttt{OMP\_PROC\_BIND=close} and \texttt{OMP\_PLACES=cores}, and only the GEMV size is
increased to $n{=}4000$.

%\paragraph{How speed is measured (GPU).}
CUDA is measured using \texttt{cudaEvent}, where the kernel is launched 20 times in a row for
AXPY and five times for GEMV and GEMM, and the mean elapsed time was taken as the time per
call.  The sweep sizes are AXPY $m\in\{2^{20},2^{22},2^{24}\}$,
GEMV $n\in\{1024,4096,8192\}$ and GEMM $n\in\{512,1024\}$, and Table~\ref{tab:cuda}
lists the largest kernel sizes.
The relative error on the GPU side has been changed to a scheme that distills the
difference as an expansion in the component type. The conventional scheme of
reconstructing into a \texttt{long double} and subtracting depends on its width
(64 bits on x86, 113 bits on aarch64), which would make the accuracy columns of H100 and
GB10 incomparable.

%\paragraph{Operation counts and the generation of the tables.}
%The operation counts of Table~\ref{tab:dsqflop} are not theoretical values but
%the actual numbers of \texttt{add/sub/mul/neg/div/sqrt/fma} counted at run time
%by a counting backend.

\subsection{How the test problems are given}\label{sec:testdata}
We describe the concrete numbers given to each kernel.
Below, $\rho$ denotes a uniform random number on the interval $(-1,1)$, $\mu$ a uniform
random number on $[1,2)$, and $b$ the number of explicitly stored fraction bits of the underlying
floating-point type ($b=52$ for double and $b=23$ for float).
All random numbers come from \texttt{rand()} with a fixed seed
(\texttt{20260709} for AXPY/GEMV/GEMM, \texttt{20260710} for division and square root and
\texttt{20260803} for CUDA), so the same binary produces the same inputs and the same
results however many times it is run.

Multiple-precision random numbers were generated as follows:
For the accuracy measurements, the value of one logical element was first formed on MPFR
with 600 bits as
\[
  v=\rho_0+\rho_1 2^{-b}+\rho_2 2^{-2b}+\rho_3 2^{-3b}
\]
and rounded into $K$ components to yield a non-overlapping multiple-precision random
value. That is, $|v|\lesssim1$ and each component carries the weight
$1,2^{-b},2^{-2b},\dots$, and the components are filled with random bits down to the
last word (so that no component becomes $0$ and misses an opportunity for cancellation).
For the timing measurements, the components were formed directly without MPFR as
$v_0=\rho$ and $v_{c+1}=v_c\cdot\rho\cdot2^{-b}$.
Since $2^{-b}=2u$, where $u=2^{-p}$ is the unit roundoff of the precision
$p=b+1$ including the hidden bit, the timing data generated in this direct way do not necessarily
satisfy the non-overlapping condition~\eqref{eq:nonoverlap} assumed by the error
analysis. This does not affect the speed comparison because these data are used only
for timing and not for the accuracy evaluation; all accuracy measurements use the
MPFR-based non-overlapping generation above.
None of the four variants compared has a value-dependent branch (except TW/QW of the
existing Q), so the execution time does not depend on the input values and this
difference does not affect the speed ratios.

The operands $a$ and $b$ of the division $a/b$ and the square root $\sqrt{a}$ are
generated independently as follows:
\[\begin{split}
  a&=\Bigl|\,\mu\,2^{e}+\sum_{j=1}^{3}\rho_j 2^{\,e-jb}\,\Bigr|,\\
  &\qquad e\in\{-8,\dots,8\}\ \text{(uniform integer)} .
\end{split}\]
Thus, $a$ and $b$ are both positive with $|a|,|b|\in[2^{-8},2^{9})$, their exponents
differ by at most $17$ bits, and they are non-overlapping $K$-component values filled
with random bits down to the last word.
The square root uses the same array $a$ (it is always defined because $a>0$).
Therefore, division and square root are compared on the same operands.
The accuracy measurement uses $n=8192$ elements and the timing measurement
$n=131072$ elements; on the timing side, the components of $a$ are $\mu_1 2^{-tb}$ and
those of $b$ are $\mu_2 2^{-tb}$ ($t=0,\dots,K-1$), i.e.,\ simple values with all
components positive.

In the computation of AXPY $\mathbf{y}\leftarrow\alpha \mathbf{x}+\mathbf{y}$, $\alpha$ is
a single multiple-precision scalar: the multiple-precision random value described above for the
accuracy measurement, and the fixed value $\alpha_t=0.9\cdot2^{-tb}$ ($t=0,\dots,K-1$)
for the timing measurement (on CUDA, $\alpha_0=0.9$ and $\rho\,2^{-tb}$ below it).
Care is needed to match the spacing of the components to the component type: reusing the
$2^{-52}$ spacing of the double-precision base for the single-precision base would make
the lowest component fall below the smallest normal number of float,
$1.18\times10^{-38}$, so that $\alpha$ itself would be subnormal.
$\mathbf{x}$ and $\mathbf{y}$ are real vectors of length $n$ whose elements are all the
multiple-precision random values above.
$\mathbf{y}$ is updated in place starting from the same initial copy for each variant, so
that the four variants see exactly the same input.

In the computation of GEMV $\mathbf{y}\leftarrow A\mathbf{x}+\mathbf{y}$, $A$ is an
$n\times n$ matrix, and $\mathbf{x}$ is a real vector of length $n$, both of whose elements
are multiple-precision random values.
$A$ is stored as a column-major SoA split by component.
Since the accumulator starts from $\mathbf{y}=0$ and overwrites $\mathbf{y}$, the initial
value of $\mathbf{y}$ does not affect the result.% ($\beta=0$).
Only the CUDA version stores $A$ row-major, because one warp handles one row; there we
compute $\mathbf{y}\leftarrow A\mathbf{x}+\mathbf{y}$ (the added $\mathbf{y}$ is a
multiple-precision random value as well) and fold the partial sums of the 32 lanes of a
warp with the same MAC using an exact multiplication by $1$.

In the computation of GEMM $C\leftarrow AB+C$, $A$ and $B$ are both $n\times n$
matrices of multiple-precision random values, and since we start from $C = 0$ the
computation is effectively $C=AB$ (no matrix is added).
The inner loop is exactly the MAC $c_{ij}\leftarrow a_{ik}b_{kj}+c_{ij}$, giving a
dependency chain of length $n$ per output element.

When an ill-conditioned inner product is used in a cancellation test, to determine the
effect of the error bound of the proposed FMA being relative to $(|xy|+|c|)$, we build an
inner product whose result cancels severely, in the shape of a GEMV.
We set $n=512$ and pair the adjacent elements of $x$ as
\[\begin{split}
  x_{2k}&=\xi_k,\qquad x_{2k+1}=-\xi_k(1+\delta),\\%\qquad
  \xi_k&=(1+0.5\rho)\,2^{s_k},\quad s_k\in[-40,40]
\end{split}\]
and fill every row of $A$ with the same vector $a_j=1.5+0.4\rho$ (where $a_{2k}=a_{2k+1}$).
Then, $y_i$ is the same inner-product value for every $i$, and only $-\delta\,a_{2k}\xi_k$
remains as the contribution of each pair; therefore, the condition number is
$\sum_j|a_jx_j|/|\sum_j a_jx_j|\approx2/\delta$.
The values $\delta=10^{-3},10^{-8},10^{-13}$ correspond to
$2\times10^{3},2\times10^{8},2\times10^{13}$, respectively.
Only in this test did we place a value in the high-order word alone and set the low-order
words to $0$; the input itself is exactly representable in the component type, and the
cancellation occurs only during accumulation.

The ``compute-only MAC chain'' is a chain of $20{,}000$ repetitions of
$c\leftarrow x\cdot y+c$ with the multiple-precision random values $x,y$ and the initial
value $c$ kept in the registers to remove the influence of the bandwidth limit (the
only memory traffic is at the very beginning and the very end).
For each thread $i$, we add $(i\bmod 8)\cdot10^{-6}$ to $x_0$, which prevents the compiler
from eliminating the loop by constant folding.

%--------------------------------------%
% Benchmark tests
%--------------------------------------%
\subsection{Acceleration of division and square root}\label{sec:divsqrt}

%\subsubsection{Motivation and algorithms}
In multiple-precision division and square root, the inner loop consists entirely of MACs
of the form $x\cdot y+c$.
Therefore, the proposed FMA should improve the performance not only of a single MAC or of
AXPY/GEMV/GEMM but also of division and square root.
In this section, we verify this for all six types DS/TS/QS/DD/TD/QD.

\begin{description}
\item[Division (long division)]
We take $K+1$ digits of the quotient $q = a/b$ and correct the residual at each digit:
\[\begin{split}
r &\leftarrow a, \\%\qquad
q_i &\leftarrow \frac{r_0}{b_0},\qquad
r \leftarrow r - q_i\cdot b \quad (i=0,\dots,K-1)\\ %\qquad
q_K&\leftarrow \frac{r_0}{b_0},
\end{split}\]
and finally distill the $K+1$ values $q_0,\dots,q_K$ into $K$ components,
($K$ bubble sweeps).
Note that the correction $r\leftarrow r-q_i b$ is exactly a MAC and the multiplier
$q_i$ is a constant.

\item[Square root (Newton iteration for $1/\sqrt{a}$)]
\[\begin{split}
r_0&\leftarrow \frac{1}{\sqrt{a_0}},\qquad h\leftarrow \frac{a}{2},\\%\qquad
e &\leftarrow \frac{1}{2} - h\cdot r^2,\qquad
r \leftarrow r + r\cdot e,\\ %\qquad
\sqrt{a}&\simeq a\cdot r .
\end{split}
\]
The updates of both $e$ and $r$ are MACs. The iteration starts from about $p-1$
correct bits and doubles them each time. Therefore, for the target of $K\cdot p$ bits, we use two
iterations for $K=2,3$ and three for $K=4$ (with only two iterations for $K=4$, one would obtain
$208<212$ bits for QD and $92<96$ bits for QS, which is insufficient).
The iteration count is common to the three methods, and therefore does not affect the
speed ratios.
\end{description}

%\paragraph{The three methods}
We implemented and compared the MAC in these two algorithms in the following three ways:
\begin{itemize}
\item \textbf{BF (baseline)}: the two stages of branch-free \texttt{mul}$+$\texttt{add}
      (Zhang--Aiken).  The product is rounded back into $K$ components and then added.
\item \textbf{FMA (proposed)}: the proposed branch-free FMA that merges the distillation
      into one stage.
\item \textbf{Exact FMA}: the result-relative distillation FMA.
\end{itemize}
The squaring $r\cdot r$ and the final $a\cdot r$ use the same \texttt{mul} in all three
methods (a real library always has a cheap \texttt{mul}).
The difference appears only in the MAC part, which is the contribution of
the proposed FMA.

%\paragraph{The zero-folded versions \texttt{mul\_d} / \texttt{fma\_d}}
Because the multiplier $q_i$ in the correction stage of the division is a scalar, we use
the functions \texttt{mul\_d} / \texttt{fma\_d} in which $y_1=\dots=y_{K-1}=0$ has been
folded in.
Because $\mathrm{TwoProd}(x,0)=(0,0)$ and $\mathrm{TwoSum}(a,0)=(a,0)$ are both
exact, the folding did not change rounding.
\texttt{fma\_d}$(x,q,c)$ agrees bitwise with \texttt{fma}$(x,\{q,0,\dots\},c)$
(six types $\times$ 20{,}000 trials with mismatch$=0$).
The BF side uses \texttt{mul\_d} in the same manner, making the comparison fair. This also
matches the configuration used by \texttt{\_bncneon\_rqd\_div} in BNCmatmul and
related routines (\texttt{mul\_d}$+$\texttt{sub}).

%\subsubsection{Implementation}
Division and square root are written as backend-independent C++ templates and instantiated
from the same source for scalar execution, NEON (\texttt{float64x2\_t},
\texttt{float32x4\_t}), SVE2 (fixed-length \texttt{svfloat64\_t},
\texttt{svfloat32\_t}), AVX2, AVX-512, and CUDA (scalar per thread).
Since all of them reduce to the same IEEE-754 basic operations and a genuine hardware
FMA, they are designed to produce identical rounding behavior across backends.
We confirmed by measurement that the NEON / SVE2 / CUDA (device) kernels agree bitwise
with the scalar execution of the same source ($4096$--$8192$ points each, mismatch$=0$).
The relative-error tables of NEON and SVE2 agree in all 36 rows as well.
\texttt{-ffp-contract=off} on the host and \texttt{--fmad=false} on nvcc are mandatory.

%\paragraph{Context dependence of the normalization gates and the div/sqrt-safe versions (machine-proved).}
A gate can be changed in two directions.  Turning an unproved \FTS{} into \TS{} (safe)
makes it provable without any loss of accuracy, because \TS{} is unconditionally exact
($+3$ flops per gate), whereas turning a \TS{} into an \FTS{} (dangerous) is exact only
in a context where the premise holds.
The standard MAC can prove the \FTS{} premise of $(B,m_1)$/$(B,A_1)$, but the FMA inside
division and square root receives intermediate values whose non-overlapping property is
not guaranteed, such as the residual of the iteration, so the premise fails and applying
this optimization degrades the relative error.
Since the non-overlapping property of the inputs cannot be assumed in div/sqrt, we
therefore use a ``safe version'' in which every \FTS{} is replaced by \TS{}
(DW $20$ / TW $84$ / QW $206$ flops; the number of normalization passes is the same as
in the standard version).  This is why the minimal flop counts of the proved version
differ.

%\subsubsection{Dynamic operation counts}
Table~\ref{tab:dsqflop} lists the actual number of operations executed (counting each of
add/sub/mul/neg/div/sqrt/fma as one).
They were counted at runtime with a counting backend and are the exact values of the
code for the benchmark runs.
They do not depend on the component type (float/double) and are determined only by
the number of components $K$.

\begin{table}[h]\centering
\caption{Dynamic operation count per element for division and square root (Exact is the fully distilled version, $n{=}K{+}1$)}\label{tab:dsqflop}
\setlength{\tabcolsep}{3pt}
\small
\begin{tabular}{lrrrc}
\toprule
& \multicolumn{4}{c}{division} \\%& \multicolumn{4}{c}{square root}\\
\cmidrule(lr){2-5} \\
$K$ & BF & FMA & Exact & \textbf{BF/FMA} \\
\midrule
2 (DS,DD) &  87 &  67 &  425 & \textbf{1.30} \\
3 (TS,TD) & 319 & 277 & 1744 & \textbf{1.15} \\
4 (QS,QD) & 785 & 749 & 4881 & \textbf{1.05} \\
\midrule
& \multicolumn{4}{c}{square root}\\
\cmidrule(lr){2-5} \\
2 (DS,DD) &  160 &  120 &  828 & \textbf{1.33}\\
3 (TS,TD) &  536 &  482 & 2378 & \textbf{1.11}\\
4 (QS,QD) & 1760 & 1730 & 7622 & \textbf{1.02}\\
\bottomrule
\end{tabular}
\end{table}

%\subsubsection{Speed}
Table~\ref{tab:dsqspeed} lists the speed ratio of the proposed FMA to the BF
baseline. The CPU used one core (\texttt{taskset}) with $n=131072$ elements, and the GPU
$n=4{,}194{,}304$ elements.
The CPU side used four backends: NEON and SVE2 on Arm, and AVX2 and AVX-512 on x86. The
GPU side used two machines: NVIDIA GB10 (sm\_121, 48 SM, LPDDR5X) and NVIDIA H100 NVL
(sm\_90, 132 SM, HBM3). Both were evaluated from the same source using the same arguments.
The backend is a thin layer that maps the six operations---addition, subtraction,
multiplication, division, square root, and hardware FMA---one-to-one onto the intrinsics of
each ISA, and the body of the division and square-root algorithms is completely common to
the six backends.

\begin{table*}[h]\centering
\caption{Speed ratio FMA/BF of division and square root ($>1$ means faster)}\label{tab:dsqspeed}
\setlength{\tabcolsep}{3pt}
\small
\begin{tabular}{llrrrrrr}
\toprule
operation & environment & DS & TS & QS & DD & TD & QD\\
\midrule
\multirow{8}{*}{division}
 & NEON              & \textbf{1.50} & 1.37 & 1.17 & \textbf{1.47} & 1.36 & 1.17\\
 & SVE2 (VL$=$128)   & \textbf{1.56} & 1.38 & 1.19 & \textbf{1.54} & 1.38 & 1.19\\
 & AVX2              & 1.61 & 1.28 & \textbf{1.83} & \textbf{1.61} & 1.28 & 1.10\\
 & AVX-512           & 1.48 & 1.35 & \textbf{1.84} & \textbf{1.48} & 1.34 & 1.15\\
 & GB10 element-wise & 1.00 & 1.00 & 1.00 & \textbf{1.29} & 1.13 & 1.04\\
 & GB10 compute-only & \textbf{1.24} & 1.13 & 1.04 & \textbf{1.29} & 1.13 & 1.04\\
 & H100 element-wise & \textbf{1.08} & 1.08 & 1.05 & 1.02 & \textbf{1.03} & 1.03\\
 & H100 compute-only & \textbf{1.25} & 1.11 & 1.04 & \textbf{1.24} & 1.13 & 1.04\\
\midrule
\multirow{8}{*}{square root}
 & NEON              & \textbf{1.53} & 1.30 & 1.07 & \textbf{1.54} & 1.30 & 1.07\\
 & SVE2 (VL$=$128)   & \textbf{1.59} & 1.29 & 1.13 & \textbf{1.57} & 1.29 & 1.13\\
 & AVX2              & \textbf{1.47} & 1.20 & 1.17 & \textbf{1.46} & 1.21 & 1.03\\
 & AVX-512           & \textbf{1.49} & 1.26 & 1.21 & \textbf{1.44} & 1.27 & 1.10\\
 & GB10 element-wise & 1.00 & 1.00 & 1.01 & \textbf{1.25} & 1.09 & 1.00\\
 & GB10 compute-only & \textbf{1.25} & 1.09 & 1.03 & \textbf{1.25} & 1.09 & 1.00\\
 & H100 element-wise & \textbf{1.14} & 1.08 & 1.02 & \textbf{1.11} & 1.07 & 1.01\\
 & H100 compute-only & \textbf{1.22} & 1.09 & 1.01 & \textbf{1.25} & 1.08 & 1.00\\
\midrule
\multicolumn{2}{l}{ref.: operation-count ratio BF/FMA (division)} & 1.30 & 1.15 & 1.05 & 1.30 & 1.15 & 1.05\\
\multicolumn{2}{l}{ref.: operation-count ratio BF/FMA (square root)} & 1.33 & 1.11 & 1.02 & 1.33 & 1.11 & 1.02\\
\bottomrule
\end{tabular}
\end{table*}

\begin{table*}[h]\centering
\caption{Relative error of division and square root (against MPFR 600-bit, max / mean, $n=8192$).
Bold marks the best value in each row.
The values are those of the NEON implementation and agree with the SVE2 implementation
in all 36 rows (the CUDA implementation uses the same algorithm but slightly different
numbers because its test inputs differ).}\label{tab:dsqacc}
\setlength{\tabcolsep}{3pt}
\small
\begin{tabular}{llll}
\toprule
type & BF (mul$+$add) & FMA (proposed) & Exact FMA\\
\midrule
\multicolumn{4}{l}{\textbf{Division} $c=a/b$}\\
DS & \ee{8.14}{-15} / \ee{9.15}{-16} & \ee{9.27}{-15} / \ee{1.22}{-15} & \textbf{\ee{1.76}{-15}} / \ee{4.37}{-16}\\
TS & \ee{3.43}{-22} / \ee{2.47}{-23} & \ee{5.46}{-22} / \ee{3.78}{-23} & \textbf{\ee{5.04}{-23}} / \ee{8.48}{-24}\\
QS & \ee{1.65}{-29} / \ee{6.59}{-31} & \ee{2.31}{-29} / \ee{1.11}{-30} & \textbf{\ee{1.53}{-30}} / \ee{1.67}{-31}\\
DD & \ee{2.24}{-32} / \ee{3.11}{-33} & \ee{3.39}{-32} / \ee{4.14}{-33} & \textbf{\ee{6.04}{-33}} / \ee{1.49}{-33}\\
TD & \ee{2.15}{-48} / \ee{1.64}{-49} & \ee{6.65}{-48} / \ee{2.48}{-49} & \textbf{\ee{3.34}{-49}} / \ee{5.39}{-50}\\
QD & \ee{1.93}{-64} / \ee{8.11}{-66} & \ee{2.35}{-64} / \ee{1.27}{-65} & \textbf{\ee{1.89}{-65}} / \ee{2.04}{-66}\\
\midrule
\multicolumn{4}{l}{\textbf{Square root} $c=\sqrt{a}$}\\
DS & \ee{1.37}{-14} / \ee{1.70}{-15} & \ee{1.37}{-14} / \ee{1.70}{-15} & \textbf{\ee{1.02}{-14}} / \ee{1.57}{-15}\\
TS & \ee{7.06}{-22} / \ee{6.57}{-23} & \ee{7.06}{-22} / \ee{6.69}{-23} & \ee{7.06}{-22} / \textbf{\ee{5.92}{-23}}\\
QS & \ee{3.17}{-29} / \ee{2.48}{-30} & \ee{3.17}{-29} / \ee{2.42}{-30} & \ee{3.17}{-29} / \textbf{\ee{2.32}{-30}}\\
DD & \ee{3.87}{-32} / \ee{5.71}{-33} & \ee{3.87}{-32} / \ee{5.71}{-33} & \ee{3.87}{-32} / \textbf{\ee{5.33}{-33}}\\
TD & \textbf{\ee{3.75}{-48}} / \ee{4.26}{-49} & \ee{3.94}{-48} / \ee{4.38}{-49} & \ee{3.79}{-48} / \textbf{\ee{3.87}{-49}}\\
QD & \textbf{\ee{3.85}{-64}} / \ee{2.95}{-65} & \ee{4.03}{-64} / \ee{2.99}{-65} & \textbf{\ee{3.85}{-64}} / \textbf{\ee{2.75}{-65}}\\
\bottomrule
\end{tabular}

\end{table*}

The CPU exceeded the operation-count ratio for the following reasons.

The measured division for $K{=}2$ is DS/DD \spd{1.48}--\spd{1.61}, above the
operation-count ratio $1.30$.
On the BF path two sequential \FTS{} chains---the normalization of
\texttt{mul\_d} and that of \texttt{add}---are placed in series, whereas the proposed FMA
folds them into one.
In addition to the reduction of the operation count, the dependency chain becomes
shorter and the ILP improves.
For $K\ge3$, on the other hand, the div/sqrt-safe version cannot use any \FTS{} and
keeps the same number of normalization passes as the standard version, so the
operation-count ratio itself rapidly approaches 1 (division $1.15/1.05$, square root
$1.11/1.02$) and the measured level also drops to TS/TD \spd{1.20}--\spd{1.35} and
QS/QD \spd{1.03}--\spd{1.84}.
\textbf{The advantage is most consistent for DW (DS/DD); for $K\ge3$ it is generally
smaller or backend-dependent, although QS division on x86 is a notable exception.}

%\paragraph{x86 shows the same trend as Arm.}
The x86 AVX2 (\texttt{\_\_m256d} 4 / \texttt{\_\_m256} 8) and AVX-512 (8 / 16) show the
same trend: division \spd{1.10}--\spd{1.84} and square root \spd{1.03}--\spd{1.49}
(Table~\ref{tab:dsqspeed}).
The broad trend is therefore similar across the Arm and x86 ISAs, although the
magnitude of the speed-up remains backend-dependent.
The backend merely maps the six operations---addition, subtraction, multiplication,
division, square root and hardware FMA---one-to-one onto the x86 intrinsics, and the
algorithm is common to the Arm version.
Indeed the relative errors agree bitwise with NEON in all 36 values, and the comparison
of the SIMD FMA against the scalar reference showed no mismatch in any of the
6 types $\times$ \{division, square root\}.
Note that AVX-512 does not consistently beat AVX2 (DS division is \spd{1.48} against
\spd{1.61}, below it):
doubling the lane width does not double the throughput of division and square root
themselves, so the effect of the width does not appear as directly as in the MAC-dominated
BLAS kernels (Table~\ref{tab:time_x86}).

%\paragraph{H100: all six types move to the bandwidth-bound side}
Running the same source with the same arguments on an H100 NVL (sm\_90, 132 SM, HBM3),
the FMA/BF of the element-wise kernels becomes \spd{1.02}--\spd{1.08} for division and
\spd{1.01}--\spd{1.14} for square root, i.e.\ \textbf{almost 1}.
The FP64 of the H100 is as strong as $1/2$ of its FP32, so the memory side saturates
first and all six types fall on the bandwidth-bound side.
The compute-only chain, on the other hand, gives \spd{1.04}--\spd{1.25} for division and
\spd{1.00}--\spd{1.25} for square root, \textbf{reproducing well the operation-count
ratios (division $1.30/1.15/1.05$, square root $1.33/1.11/1.02$)}.
In other words, the advantage in operation count appears exactly as the operation-count
ratio, independently of the GPU generation and the memory type, and only the degree to
which it appears in the wall-clock time is decided by the compute/bandwidth balance of
the kernel.
Since the operation-count ratio of the div/sqrt-safe version itself approaches 1 with
growing $K$, however, \textbf{the advantage on the GPU also remains only for DW}
(QS/QD are at most \spd{1.05} both element-wise and compute-only).

%\paragraph{Why the three single-precision-based types are \spd{1.00} on GB10}
The reason why the three single-precision-based types are \spd{1.00} on GB10 can be
explained as follows.

The element-wise kernel only reads and writes $3K$ component arrays, so its arithmetic
intensity is low.  The measured effective bandwidth is DS $235$ / TS $232$ / QS $228$ GB/s,
reaching $84$--$86\%$ of the LPDDR5X peak of GB10 (about $273$ GB/s), i.e.\ it is
bandwidth-bound.  The three double-precision-based types become compute-bound because the
FP64 of GB10 is weak (DD $131$ / TD $52$ / QD $26$ GB/s), so the reduction of the
operation count appears directly in the wall-clock time.
This is the same picture as the GEMM of Section~\ref{sec:cuda}.
If only the arithmetic is extracted with a register-resident compute-only chain
(\texttt{chdiv} / \texttt{chsqrt}), which does not touch the global memory inside the loop,
all six types agree with the theoretical operation-count ratio, as shown in
Table~\ref{tab:dsqspeed}.
Thus, the advantage in the operation count is real on the GPU for division and square root, as
well, but since the operation-count ratio of the div/sqrt-safe version approaches 1 with
growing $K$, its absolute magnitude is largest for DW.

%\paragraph{Exact FMA}
The Exact/BF of the fully distilled version ($n{=}K{+}1$) is, for division,
NEON \spd{0.25}--\spd{0.29}, SVE2 \spd{0.11}--\spd{0.27},
AVX2 \spd{0.11}--\spd{0.35}, AVX-512 \spd{0.13}--\spd{0.32} and
H100 \spd{0.18}--\spd{0.56}; for square root,
NEON \spd{0.27}--\spd{0.34}, SVE2 \spd{0.12}--\spd{0.27},
AVX2 \spd{0.06}--\spd{0.28}, AVX-512 \spd{0.06}--\spd{0.36} and
H100 \spd{0.22}--\spd{0.34}
(corresponding to the growth in the dynamic operation count of the Exact caused by the one
additional sweep: DIV $1384\to1744$ (TW) / $4041\to4881$ (QW),
SQRT $1898\to2378$ / $6362\to7622$).
The decrease in SVE2 was the largest because the 36 fixed-length vectors (QD) required by the
distillation sweeps did not fit into the registers and spills, which is the same phenomenon as that for
MAC, as described in Section~\ref{sec:sve2}. The proposed FMA and BF have no expansion arrays and
are therefore unaffected.

%\subsubsection{Relative error}
Table~\ref{tab:dsqacc} lists the errors relative to the MPFR 600-bit reference data
($n=8192$, maximum/mean).
Here, $\varepsilon_K=2^{-Kp}$ is the unit roundoff for each type
(DS \ee{3.55}{-15}, TS \ee{2.12}{-22}, QS \ee{1.26}{-29}, DD \ee{1.23}{-32},
TD \ee{1.37}{-48}, QD \ee{1.52}{-64}).

%\paragraph{Division}
Regarding the relative error of division measured in units of $\varepsilon_K$, the
maximum error is $1.3$--$2.3$ for BF and $1.5$--$4.9$ for the proposed FMA.
The proposed FMA has a slightly larger error than BF (by a factor $1$--$2$, at most $3$
for TD). This is because BF is normalized once when the product is rounded back into
$K$ components in the correction stage, whereas the proposed FMA discards more of the
low-order cross terms.
Both, however, stay within a few times $\varepsilon_K$; for these test problems, their
accuracy is therefore practically comparable.
The fully distilled Exact FMA attains a maximum error of $0.12$--$0.50\,\varepsilon_K$
for all six types (DS/DD $0.49$, TS/TD $0.24$, QS/QD $0.12$), achieving in division the
same ``error level corresponding to correct rounding'' (measured within the tested
range) as for the single operation of Section~\ref{sec:exact}.

%\paragraph{Square root}
As for the relative error of the square root, the errors of the three methods are
almost identical (for DS and DD the reported maximum errors coincide).
Since the Newton iteration converges well below $\varepsilon_K$, the final error is
determined solely by the rounding of the last multiplication $a\cdot r$, which is common
to the three methods ($2.5$--$3.9$ times $\varepsilon_K$).
For the square root, therefore, the proposed FMA is (for DW) \spd{1.44}--\spd{1.59}
faster with almost no loss of accuracy---an almost pure gain.
For the same reason there is little point in building a square root with the Exact FMA
(the final multiplication dominates the error).

%\subsubsection{Summary}

These results show that the proposed FMA can accelerate not only the MAC-type kernels
(AXPY/GEMV/GEMM) but also division and square root.
Division is \spd{1.47}--\spd{1.61} faster on the CPU (DW) and \spd{1.24}--\spd{1.29}
faster in the compute-only chains on the GPU with accuracy equivalent to BF, and the
square root is \spd{1.44}--\spd{1.59} faster on the CPU (DW) with accuracy practically
identical to BF.
Unlike the MAC family, however, the div/sqrt-safe version cannot use any \FTS{} and
keeps the same number of normalization passes as the standard version, so the
operation-count ratio rapidly approaches 1 with $K$ (division $1.30/1.15/1.05$, square
root $1.33/1.11/1.02$).
\textbf{The advantage in division and square root is therefore most consistent for DW
(DS/DD); for $K\ge3$ it is generally smaller or backend-dependent.  This differs from the
behavior of the MAC kernels.}

%--------------------------------------%
% Basic linear algebra
%--------------------------------------%
\subsection{Basic linear algebra: AXPY, GEMV, and GEMM}\label{sec:blas}

Here, we evaluate the accuracy and speed of the three kernels
AXPY ($\mathbf{y}\leftarrow\alpha \mathbf{x}+\mathbf{y}$, $n=10^6$), GEMV ($n=2048$), and GEMM ($n=512$),
focusing on the result of applying the proposed FMA.

\subsubsection{Accuracy (relative error against MPFR)}

Table~\ref{tab:acc} lists the maximum and mean relative errors for all six types
(DS/TS/QS/DD/TD/QD) in a single table.
For the backend/type combinations subjected to the cross-backend check, BF, the proposed
FMA, and the Exact FMA agree bitwise; only the Q column can differ because it uses the
default implementation of each ISA.  Table~\ref{tab:acc} therefore summarizes the
accuracy results, with the exact backend coverage stated in its caption.
All six types maintained their working accuracy
(DS$\sim$14, TS$\sim$21, QS$\sim$28, DD$\sim$31, TD$\sim$46, QD$\sim$61 digits).
The proposed FMA is comparable to existing methods (Q/BF); in the worst case, it can
degrade by approximately one digit (\spd{8.0} against BF for QS GEMM), originating from the
$(|xy|+|c|)$-relative semantics; however, in the mean error, it remains at the same level as the BF
(it can be substituted without a loss of accuracy).
The fully distilled Exact FMA was the best in 15 of the 18 rows for the maximum relative
error and 17 of 18 for the mean relative error, and is overwhelming in AXPY, where the
rounding of a single operation dominated this error
(DS \ee{3.1}{-12}$\to$\ee{1.7}{-15}, QS \ee{4.7}{-27}$\to$\ee{1.5}{-30},
TD \ee{2.3}{-46}$\to$\ee{3.3}{-49}, QD \ee{4.1}{-61}$\to$\ee{1.8}{-65}).
The rows in which it is not the best (QS GEMM, DD GEMV and TD GEMM for the maximum error,
TD GEMM for the mean error) correspond to cases where the maximum error is decided not by the
rounding of each fused operation, but rather by the cancellation in the accumulation along $k$,
which the result-relative semantics does not improve either.

\begin{table*}[htb]\centering
\caption{Summary of the relative errors of all six types (against MPFR 600-bit).
Bold marks the best value in each row.
The upper block is the single-precision base DS/TS/QS and the lower the
double-precision base DD/TD/QD.
The values of BF, FMA, and Exact agree bitwise across every backend that implements the
type---AVX2 / AVX-512 / NEON / SVE2 / CUDA for the double-precision base and
AVX2 / AVX-512 / NEON / CUDA for the single-precision base.
Only the Q column is the default implementation of the existing library for each ISA and
does not agree (the values are those of NEON; refer to the text).}\label{tab:acc}
\footnotesize
\setlength{\tabcolsep}{3pt}
\begin{tabular}{ll cccc cccc}
\toprule
 & & \multicolumn{4}{c}{maximum relative error} & \multicolumn{4}{c}{mean relative error}\\
\cmidrule(lr){3-6}\cmidrule(lr){7-10}
type & kernel & Q & BF & FMA & Exact & Q & BF & FMA & Exact\\
\midrule
DS & AXPY & \ee{4.1}{-12} & \ee{3.6}{-12} & \ee{3.1}{-12} & \textbf{\ee{1.7}{-15}} & \ee{6.1}{-15} & \ee{5.1}{-15} & \ee{5.7}{-15} & \textbf{\ee{4.3}{-16}}\\
DS & GEMV & \ee{1.3}{-12} & \ee{1.2}{-12} & \ee{1.5}{-12} & \textbf{\ee{5.6}{-13}} & \ee{3.9}{-14} & \ee{3.4}{-14} & \ee{4.3}{-14} & \textbf{\ee{1.7}{-14}}\\
DS & GEMM & \ee{2.0}{-10} & \ee{8.9}{-11} & \ee{1.4}{-10} & \textbf{\ee{5.1}{-11}} & \ee{4.7}{-14} & \ee{4.0}{-14} & \ee{5.0}{-14} & \textbf{\ee{2.6}{-14}}\\
TS & AXPY & \ee{6.1}{-21} & \ee{3.0}{-20} & \ee{5.9}{-20} & \textbf{\ee{5.0}{-23}} & \ee{2.0}{-23} & \ee{4.7}{-23} & \ee{8.0}{-23} & \textbf{\ee{8.6}{-24}}\\
TS & GEMV & \ee{1.1}{-19} & \ee{3.2}{-19} & \ee{2.1}{-19} & \textbf{\ee{4.6}{-20}} & \ee{7.7}{-22} & \ee{2.4}{-21} & \ee{2.3}{-21} & \textbf{\ee{6.9}{-22}}\\
TS & GEMM & \ee{2.5}{-19} & \ee{6.6}{-19} & \ee{9.3}{-19} & \textbf{\ee{2.2}{-19}} & \ee{3.6}{-22} & \ee{7.2}{-22} & \ee{9.5}{-22} & \textbf{\ee{3.1}{-22}}\\
QS & AXPY & \ee{4.1}{-27} & \ee{1.4}{-26} & \ee{4.7}{-27} & \textbf{\ee{1.5}{-30}} & \ee{8.8}{-30} & \ee{8.7}{-30} & \ee{6.8}{-30} & \textbf{\ee{1.7}{-31}}\\
QS & GEMV & \ee{4.5}{-27} & \ee{3.9}{-27} & \ee{3.7}{-27} & \textbf{\ee{1.4}{-27}} & \ee{3.5}{-29} & \ee{3.8}{-29} & \ee{4.7}{-29} & \textbf{\ee{1.1}{-29}}\\
QS & GEMM & \ee{2.0}{-25} & \textbf{\ee{5.5}{-26}} & \ee{4.4}{-25} & \ee{1.9}{-25} & \ee{4.2}{-29} & \ee{3.5}{-29} & \ee{7.1}{-29} & \textbf{\ee{1.9}{-29}}\\
\midrule
DD & AXPY & \ee{2.9}{-29} & \ee{2.9}{-29} & \ee{2.1}{-29} & \textbf{\ee{6.0}{-33}} & \ee{3.1}{-32} & \ee{2.2}{-32} & \ee{2.7}{-32} & \textbf{\ee{1.5}{-33}}\\
DD & GEMV & \ee{3.3}{-30} & \textbf{\ee{2.8}{-30}} & \ee{4.7}{-30} & \ee{3.7}{-30} & \ee{9.7}{-32} & \ee{1.1}{-31} & \ee{1.5}{-31} & \textbf{\ee{6.5}{-32}}\\
DD & GEMM & \ee{4.2}{-28} & \ee{4.3}{-28} & \ee{1.1}{-27} & \textbf{\ee{3.2}{-28}} & \ee{1.7}{-31} & \ee{1.6}{-31} & \ee{2.1}{-31} & \textbf{\ee{8.0}{-32}}\\
TD & AXPY & \ee{5.3}{-47} & \ee{3.1}{-46} & \ee{2.3}{-46} & \textbf{\ee{3.3}{-49}} & \ee{1.3}{-49} & \ee{3.4}{-49} & \ee{4.3}{-49} & \textbf{\ee{5.4}{-50}}\\
TD & GEMV & \ee{1.0}{-45} & \ee{7.3}{-46} & \ee{1.1}{-45} & \textbf{\ee{3.0}{-46}} & \ee{7.2}{-48} & \ee{8.6}{-48} & \ee{1.2}{-47} & \textbf{\ee{4.2}{-48}}\\
TD & GEMM & \textbf{\ee{2.0}{-45}} & \ee{4.4}{-45} & \ee{5.8}{-45} & \ee{5.1}{-45} & \textbf{\ee{2.1}{-48}} & \ee{4.9}{-48} & \ee{6.9}{-48} & \ee{2.3}{-48}\\
QD & AXPY & \ee{5.4}{-61} & \ee{3.2}{-61} & \ee{4.1}{-61} & \textbf{\ee{1.8}{-65}} & \ee{3.4}{-64} & \ee{1.9}{-64} & \ee{2.9}{-64} & \textbf{\ee{2.0}{-66}}\\
QD & GEMV & \ee{1.9}{-62} & \ee{4.0}{-62} & \ee{2.2}{-62} & \textbf{\ee{8.1}{-63}} & \ee{3.3}{-64} & \ee{4.0}{-64} & \ee{4.5}{-64} & \textbf{\ee{1.3}{-64}}\\
QD & GEMM & \ee{6.7}{-60} & \ee{1.2}{-60} & \ee{1.6}{-60} & \textbf{\ee{8.5}{-61}} & \ee{7.9}{-64} & \ee{4.7}{-64} & \ee{5.6}{-64} & \textbf{\ee{1.6}{-64}}\\
\bottomrule
\end{tabular}
\end{table*}

%\paragraph{Cancellation (ill-conditioned inner product).}
Even at a condition number of $2\times10^{13}$, the proposed FMA degrades almost
identically to the existing ones (DD: Q \ee{1.3}{-20} $=$ FMA \ee{1.3}{-20}).
Replacing them with the proposed FMA results in no regression in accuracy.

\subsubsection{Serial performance: Arm (GB10, Cortex-X925 / NEON and SVE2)}\label{sec:neon}\label{sec:sve2}
The measurement was pinned to one performance core (Cortex-X925, 3.9\,GHz) with
\texttt{taskset} and measured on the two ISAs NEON and SVE2.
The results are listed in Table~\ref{tab:time_arm}.
The scalar EFT primitives were mapped one-to-one to
\texttt{\_bncneon\_dtwo\_sum / \_dtwo\_prod / \_dquick\_two\_sum} (which use the hardware
FMA) and plain $+,\times$ onto \texttt{vaddq\_f64 / vmulq\_f64}.
Because \texttt{float64x2\_t} of NEON holds two double-precision elements, and
\texttt{float32x4\_t} four single-precision elements, only the row-direction parallelism of
the SoA and loop stride change with type.
The compilation flags were gcc \texttt{-O3 -march=native -DBNC\_ENABLE\_NEON -funroll-loops}.
All six types (DS/TS/QS/DD/TD/QD) were implemented on both NEON and SVE2.

%\paragraph{Arm SVE2 (an inter-ISA comparison on the same machine at the same lane width).}\label{sec:sve2}
We ported the entire NEON version to Arm SVE2 (the SIMD extension of Armv9-A) and repeated
the measurement on the same machine (Cortex-X925).
The VL-agnostic primitives
\texttt{\_bncsve2\_dtwo\_sum / \_dtwo\_prod / \_dquick\_two\_sum} (with a predicate
\texttt{pg}, using the hardware FMA) provided by \texttt{include/sve2/} of BNCmatmul are
mapped one-to-one, and $+,\times$ are \texttt{svadd\_f64\_x / svmul\_f64\_x}.
The row loop is predicated with \texttt{svwhilelt\_b64} and does not require a scalar tail.
The Cortex-X925 of the measurement machine has VL$=128$ bits (\texttt{svcntd()}$=2$). Thus,
the number of double-precision lanes is two, which is the same as that in NEON.
Hence, the two column groups in Table~\ref{tab:time_arm} show an inter-ISA comparison at the
same lane width.

%\paragraph{The sizeless-type restriction and fixed-length SVE types.}
The \texttt{svfloat64\_t} of SVE is sizeless; it cannot be an array element or
structure member, and \texttt{sizeof} cannot be applied to it.
Because the Exact FMA performs distillation sweeps over an indexed expansion $v[0..35]$, it
cannot be written using a plain \texttt{svfloat64\_t}.
Therefore, our implementation builds the array with the fixed-length SVE type
\texttt{\_\_attribute\_\_((arm\_sve\_vector\_bits(N)))} of GCC/Clang
(\texttt{-msve-vector-bits=N} is mandatory) and calls the VL-agnostic primitives of the
library as they are, through the implicit conversion between the fixed-length and the
sizeless type.
The per-limb variable signature of the library (\texttt{f(pg,\&r0,\&r1,a0,a1,b0,b1)}) is
bridged with macros. At runtime, we checked whether \texttt{svcntd()} matched the number of
lanes assumed during build time.

%\paragraph{Verification of the translation.}
For the single operation $x\cdot y + c$, both the NEON implementation (all six types) and
the SVE2 implementation (DD/TD/QD) agreed bitwise with the scalar reference
($\times$ \{proposed, Exact\}, 200{,}000 trials each, no mismatch).
For the double-precision base (DD/TD/QD), these checks together with AVX2, AVX-512, and
CUDA establish bitwise agreement across all five backends.  For the single-precision base
(DS/TS/QS), the cross-backend comparison reported here covers AVX2, AVX-512, NEON, and
CUDA; the SVE2 performance results are reported separately.
The relative errors and the behavior under cancellation are otherwise consistent with
those on x86, as summarized in Table~\ref{tab:acc}.

\begin{table*}[htb]\centering
\caption{Summary of the serial measured times [s/call] on Arm (GB10, Cortex-X925 3.9\,GHz).
NEON and SVE2 (VL$=128$ bit) are shown side by side.  The number of lanes is 2 in double
precision and 4 in single precision, the same as NEON in both cases.
The upper block is the single-precision base DS/TS/QS and the lower the
double-precision base DD/TD/QD.
The fair comparison target is BF; the default implementation Q, which contains
value-dependent branches, is left out of the table and discussed in the text.
The single-precision base of SVE2 was measured with the single-precision implementation
(\texttt{\_bncsve2\_feft.h} and \texttt{\_bncsve2\_\{ds,qs\}.h}) provided by
\texttt{include/sve2/} of BNCmatmul (Section~\ref{sec:sve2}).}
\label{tab:time_arm}
\setlength{\tabcolsep}{3pt}
\footnotesize
\begin{tabular}{ll rrrr rrrr}
\toprule
 & & \multicolumn{4}{c}{\textbf{NEON}} & \multicolumn{4}{c}{\textbf{SVE2}}\\
\cmidrule(lr){3-6}\cmidrule(lr){7-10}
type & kernel & BF & FMA & Exact & FMA/BF & BF & FMA & Exact & FMA/BF\\
\midrule
DS & AXPY(1M) & \ee{6.13}{-4} & \ee{3.61}{-4} & \ee{3.24}{-3} & \spd{1.70} & \ee{6.09}{-4} & \ee{3.47}{-4} & \ee{3.49}{-3} & \spd{1.75}\\
DS & GEMV(2k) & \ee{6.16}{-3} & \ee{4.51}{-3} & \ee{1.73}{-2} & \spd{1.37} & \ee{6.20}{-3} & \ee{4.64}{-3} & \ee{1.92}{-2} & \spd{1.34}\\
DS & GEMM(512) & \ee{1.74}{-1} & \ee{1.24}{-1} & \ee{5.17}{-1} & \spd{1.40} & \ee{1.79}{-1} & \ee{1.30}{-1} & \ee{5.61}{-1} & \spd{1.37}\\
TS & AXPY(1M) & \ee{2.23}{-3} & \ee{1.24}{-3} & \ee{1.27}{-2} & \spd{1.80} & \ee{2.35}{-3} & \ee{1.23}{-3} & \ee{2.66}{-2} & \spd{1.92}\\
TS & GEMV(2k) & \ee{1.45}{-2} & \ee{1.06}{-2} & \ee{5.68}{-2} & \spd{1.36} & \ee{1.51}{-2} & \ee{1.08}{-2} & \ee{1.10}{-1} & \spd{1.40}\\
TS & GEMM(512) & \ee{3.92}{-1} & \ee{2.74}{-1} & $1.703$ & \spd{1.43} & \ee{4.19}{-1} & \ee{2.88}{-1} & $3.553$ & \spd{1.46}\\
QS & AXPY(1M) & \ee{4.55}{-3} & \ee{2.53}{-3} & \ee{3.06}{-2} & \spd{1.80} & \ee{4.93}{-3} & \ee{2.58}{-3} & \ee{7.62}{-2} & \spd{1.91}\\
QS & GEMV(2k) & \ee{2.26}{-2} & \ee{1.68}{-2} & \ee{1.26}{-1} & \spd{1.35} & \ee{2.42}{-2} & \ee{1.63}{-2} & \ee{3.07}{-1} & \spd{1.48}\\
QS & GEMM(512) & \ee{6.30}{-1} & \ee{4.31}{-1} & $4.027$ & \spd{1.46} & \ee{7.14}{-1} & \ee{4.31}{-1} & $10.19$ & \spd{1.66}\\
\midrule
DD & AXPY(1M) & \ee{1.27}{-3} & \ee{7.68}{-4} & \ee{6.43}{-3} & \spd{1.65} & \ee{1.28}{-3} & \ee{7.81}{-4} & \ee{8.26}{-3} & \spd{1.63}\\
DD & GEMV(2k) & \ee{1.24}{-2} & \ee{9.15}{-3} & \ee{3.47}{-2} & \spd{1.36} & \ee{1.25}{-2} & \ee{9.46}{-3} & \ee{4.47}{-2} & \spd{1.32}\\
DD & GEMM(512) & \ee{3.51}{-1} & \ee{2.50}{-1} & $1.048$ & \spd{1.41} & \ee{3.59}{-1} & \ee{2.62}{-1} & $1.279$ & \spd{1.37}\\
TD & AXPY(1M) & \ee{4.32}{-3} & \ee{2.38}{-3} & \ee{2.56}{-2} & \spd{1.82} & \ee{4.74}{-3} & \ee{2.56}{-3} & \ee{5.66}{-2} & \spd{1.85}\\
TD & GEMV(2k) & \ee{3.06}{-2} & \ee{2.44}{-2} & \ee{1.10}{-1} & \spd{1.26} & \ee{3.15}{-2} & \ee{2.27}{-2} & \ee{2.32}{-1} & \spd{1.39}\\
TD & GEMM(512) & \ee{7.94}{-1} & \ee{5.66}{-1} & $3.427$ & \spd{1.40} & \ee{8.44}{-1} & \ee{5.87}{-1} & $7.49$ & \spd{1.44}\\
QD & AXPY(1M) & \ee{9.00}{-3} & \ee{5.07}{-3} & \ee{6.16}{-2} & \spd{1.77} & \ee{9.87}{-3} & \ee{5.22}{-3} & \ee{1.53}{-1} & \spd{1.89}\\
QD & GEMV(2k) & \ee{4.13}{-2} & \ee{3.09}{-2} & \ee{2.50}{-1} & \spd{1.34} & \ee{4.53}{-2} & \ee{3.05}{-2} & \ee{6.27}{-1} & \spd{1.49}\\
QD & GEMM(512) & $1.302$ & \ee{8.93}{-1} & $8.069$ & \spd{1.46} & $1.454$ & \ee{8.93}{-1} & $20.85$ & \spd{1.63}\\
\bottomrule
\end{tabular}
\end{table*}

The trends in NEON and SVE2 of Arm are summarized below.

\begin{description}
\item[NEON]
The proposed FMA was \spd{1.26}--\spd{1.82} faster than the BF for all six types and all
kernels (left half of Table~\ref{tab:time_arm}).
The effect is largest for AXPY (TD \spd{1.82}, TS/QS \spd{1.80}); the GEMM gives DS \spd{1.40} / TS \spd{1.43} / QS \spd{1.46} on the single-precision
base and DD \spd{1.41} / TD \spd{1.40} / QD \spd{1.46} on the double-precision base.
Thus, on NEON the measured GEMM speed-up is nearly independent of the component count and
does not directly track the theoretical FLOP ratio.
That is, the advantage of the proposed FMA does not depend on the component type.
Against the default implementation Q, which is left out of the table, the ratio reaches
at most \spd{9.9} (TS, AXPY), but this is because the normalization renorm of Q contains
value-dependent branches; the fair comparison target remains BF
(the Q values on NEON are DS \ee{4.38}{-4} / TS \ee{1.22}{-2} / QS \ee{1.58}{-2} /
DD \ee{9.35}{-4} / TD \ee{1.02}{-2} / QD \ee{1.49}{-2} s for AXPY, and
DD \ee{2.49}{-1} / TD $1.428$ / QD $1.931$ s for GEMM(512)).
The fully distilled Exact FMA is \spd{0.15}--\spd{0.36} against the BF (corresponding to
$194/552/1178$ flops), whereas on SVE2 it drops to \spd{0.06}--\spd{0.32} because of the
stack spilling, as will be described later.

\item[SVE2]
The proposed FMA is \spd{1.32}--\spd{1.92} against BF for all six types and all kernels
(right half of Table~\ref{tab:time_arm}), the same trend as NEON.
The best GEMM values are QS \spd{1.66} on the single-precision base, followed by
QD \spd{1.63} on the double-precision base (the same combinations are \spd{1.46} on NEON).
Decomposing this difference for the QD GEMM, BF is slower on SVE2
($1.454$ s vs NEON $1.302$ s) while the proposed FMA is the same
($0.893$ s vs NEON $0.893$ s); this likely reflects differences in register availability (\texttt{z0}--\texttt{z31})
and in how effectively instruction latency is hidden.
Because the lane width is the same (2 in double, 4 in single precision), the absolute
time of the proposed FMA is roughly the same as on NEON, and for the nine
single-precision cases the difference stays within $\pm6\%$.
That is, the advantage of the proposed FMA depends neither on the ISA nor on the
component type.
The behavior under OpenMP parallelization is summarized in Section~\ref{sec:omp}.
\end{description}

%-------------------------------------------%
% Serial performance: x86
%-------------------------------------------%
\subsubsection{Serial performance: x86 (Xeon Gold 6526Y / AVX2 and AVX-512)}\label{sec:speed}
The same evaluation was conducted on x86 using AVX2 and AVX-512 (one core with
\texttt{taskset}, best of three). The results are listed in Table~\ref{tab:time_x86}.
The fair comparison target is the BF.
Both AVX2 and AVX-512 were evaluated for all six types, including the single-precision
base DS/TS/QS (the same configuration as that of the Arm/GB10 side).
All the values were measured using the fully distilled version ($n{=}K{+}1$).

%\paragraph{Four SIMD widths.}\label{sec:avx512}
Using the same Xeon Gold 6526Y, we prepared four configurations: AVX2 with 4-way parallelism
(\texttt{\_\_m256d}) and AVX-512 with 8-way (\texttt{\_\_m512d}) for the
double-precision base, and AVX2 with 8-way (\texttt{\_\_m256}) and AVX-512 with 16-way
(\texttt{\_\_m512}) for the single-precision base.
The AVX-512 version of the double-precision base is a mechanical translation of the AVX2
version: apart from \texttt{\_\_m256d}$\to$\texttt{\_\_m512d},
\texttt{\_mm256\_*}$\to$\texttt{\_mm512\_*}, the EFT primitives
\texttt{\_bncavx2\_*}$\to$\texttt{\_bncavx512\_*}, a loop stride of \texttt{i+=8} and an
alignment of 32$\to$64 bytes, and the algorithm is identical.
The two versions of the single-precision base are likewise mechanical translations of the
NEON version (\texttt{float32x4\_t}, 4-way): \texttt{vaddq\_f32}/\texttt{vmulq\_f32} are
mapped to \texttt{\_mm256\_add\_ps}/\texttt{\_mm256\_mul\_ps} (\texttt{\_mm512\_*} for
AVX-512), \texttt{\_bncneon\_f*} to \texttt{\_bncavx2\_f*} / \texttt{\_bncavx512\_f*},
and the loop stride to \texttt{i+=8} / \texttt{i+=16}.
The compilation flags are gcc \texttt{-O3 -mavx2 -mfma} (with
\texttt{-mavx512f -mavx512dq -mavx512vl -mavx512cd} added for the AVX-512 version)
\texttt{-funroll-loops}.

%\paragraph{Verification of the translation.}
For a single operation $x\cdot y+c$, both the AVX2 and AVX-512 implementations
agreed bitwise with the scalar reference for all six types
($\times$ \{proposed, Exact\}, 200{,}000 trials each, no mismatch).
Moreover, the relative errors of BF, the proposed FMA, and the Exact FMA agree completely
between the two ISAs in all 18 maximum and mean values and in all nine rows of the
cancellation test for the double-precision base, and for the single-precision base all 27
values agree between AVX2 (8-way), AVX-512 (16-way), and GB10/NEON (4-way)
(Table~\ref{tab:acc}).
This implies that the same algorithm returns the same bit pattern at 4-, 8-, and 16-way
parallelism and across Arm and x86; the rounding is invariant with respect to the SIMD
width and to the ISA.

\begin{table*}[htb]\centering
\caption{Summary of the serial measured times [s/call] on x86 (Xeon Gold 6526Y).
AVX2 (\texttt{\_\_m256d} 4 / \texttt{\_\_m256} 8) and AVX-512
(\texttt{\_\_m512d} 8 / \texttt{\_\_m512} 16) are shown side by side for all six types.
``Width'' is the AVX-512/AVX2 speed-up of the proposed FMA (the ratio of the two FMA
columns).
The fair comparison target is BF; the default implementation Q, which contains
value-dependent branches and uses different algorithms on the two ISAs, is left out of the
table and discussed in the text.}
\label{tab:time_x86}
\setlength{\tabcolsep}{2.5pt}
\footnotesize
\begin{tabular}{ll rrrr rrrr r}
\toprule
 & & \multicolumn{4}{c}{\textbf{AVX2}} & \multicolumn{4}{c}{\textbf{AVX-512}} & width\\
\cmidrule(lr){3-6}\cmidrule(lr){7-10}
type & kernel & BF & FMA & Exact & FMA/BF & BF & FMA & Exact & FMA/BF & (FMA)\\
\midrule
DS & AXPY(1M) & \ee{4.46}{-4} & \ee{4.29}{-4} & \ee{2.59}{-3} & \spd{1.04} & \ee{4.30}{-4} & \ee{4.32}{-4} & \ee{1.59}{-3} & \spd{0.99} & \spd{0.99}\\
DS & GEMV(2k) & \ee{4.02}{-3} & \ee{3.02}{-3} & \ee{1.39}{-2} & \spd{1.33} & \ee{3.07}{-3} & \ee{2.25}{-3} & \ee{9.10}{-3} & \spd{1.36} & \spd{1.34}\\
DS & GEMM(512) & \ee{9.29}{-2} & \ee{6.52}{-2} & \ee{3.70}{-1} & \spd{1.43} & \ee{7.24}{-2} & \ee{5.37}{-2} & \ee{2.33}{-1} & \spd{1.35} & \spd{1.21}\\
TS & AXPY(1M) & \ee{1.53}{-3} & \ee{9.97}{-4} & \ee{9.78}{-3} & \spd{1.54} & \ee{1.06}{-3} & \ee{6.13}{-4} & \ee{5.42}{-3} & \spd{1.72} & \spd{1.63}\\
TS & GEMV(2k) & \ee{1.30}{-2} & \ee{1.15}{-2} & \ee{4.37}{-2} & \spd{1.13} & \ee{1.06}{-2} & \ee{9.12}{-3} & \ee{2.64}{-2} & \spd{1.16} & \spd{1.26}\\
TS & GEMM(512) & \ee{2.54}{-1} & \ee{2.11}{-1} & $1.313$ & \spd{1.21} & \ee{1.75}{-1} & \ee{1.48}{-1} & \ee{7.57}{-1} & \spd{1.18} & \spd{1.42}\\
QS & AXPY(1M) & \ee{3.53}{-3} & \ee{2.65}{-3} & \ee{2.15}{-2} & \spd{1.33} & \ee{2.29}{-3} & \ee{1.71}{-3} & \ee{1.24}{-2} & \spd{1.34} & \spd{1.55}\\
QS & GEMV(2k) & \ee{2.60}{-2} & \ee{2.25}{-2} & \ee{1.02}{-1} & \spd{1.16} & \ee{2.15}{-2} & \ee{1.82}{-2} & \ee{6.43}{-2} & \spd{1.18} & \spd{1.23}\\
QS & GEMM(512) & \ee{5.03}{-1} & \ee{4.19}{-1} & $2.989$ & \spd{1.20} & \ee{3.08}{-1} & \ee{2.97}{-1} & $1.711$ & \spd{1.04} & \spd{1.41}\\
\midrule
DD & AXPY(1M) & \ee{1.01}{-3} & \ee{9.62}{-4} & \ee{5.19}{-3} & \spd{1.05} & \ee{9.81}{-4} & \ee{9.50}{-4} & \ee{3.23}{-3} & \spd{1.03} & \spd{1.01}\\
DD & GEMV(2k) & \ee{1.30}{-2} & \ee{1.02}{-2} & \ee{4.23}{-2} & \spd{1.28} & \ee{1.05}{-2} & \ee{7.97}{-3} & \ee{3.30}{-2} & \spd{1.31} & \spd{1.27}\\
DD & GEMM(512) & \ee{2.13}{-1} & \ee{1.41}{-1} & \ee{8.25}{-1} & \spd{1.51} & \ee{1.46}{-1} & \ee{1.08}{-1} & \ee{4.95}{-1} & \spd{1.35} & \spd{1.31}\\
TD & AXPY(1M) & \ee{3.17}{-3} & \ee{2.10}{-3} & \ee{1.97}{-2} & \spd{1.51} & \ee{2.20}{-3} & \ee{1.64}{-3} & \ee{1.09}{-2} & \spd{1.34} & \spd{1.28}\\
TD & GEMV(2k) & \ee{3.22}{-2} & \ee{2.89}{-2} & \ee{9.69}{-2} & \spd{1.12} & \ee{2.72}{-2} & \ee{2.38}{-2} & \ee{6.35}{-2} & \spd{1.15} & \spd{1.22}\\
TD & GEMM(512) & \ee{5.63}{-1} & \ee{4.60}{-1} & $2.682$ & \spd{1.22} & \ee{3.52}{-1} & \ee{2.97}{-1} & $1.537$ & \spd{1.18} & \spd{1.55}\\
QD & AXPY(1M) & \ee{7.09}{-3} & \ee{5.33}{-3} & \ee{4.32}{-2} & \spd{1.33} & \ee{4.65}{-3} & \ee{3.48}{-3} & \ee{2.56}{-2} & \spd{1.34} & \spd{1.53}\\
QD & GEMV(2k) & \ee{5.80}{-2} & \ee{5.28}{-2} & \ee{2.12}{-1} & \spd{1.10} & \ee{5.14}{-2} & \ee{4.27}{-2} & \ee{1.34}{-1} & \spd{1.20} & \spd{1.24}\\
QD & GEMM(512) & $1.094$ & \ee{9.04}{-1} & $6.090$ & \spd{1.21} & \ee{6.48}{-1} & \ee{6.05}{-1} & $3.453$ & \spd{1.07} & \spd{1.49}\\
\bottomrule
\end{tabular}

\smallskip
{\footnotesize The Exact FMA of the single-precision base (DS/TS/QS) was measured on a
numerical example in which all the operands were multiplied by $2^{45}$
(being a power of two, this does not move the significand by a single bit, and all 36
relative errors coincide exactly with the unscaled ones; it is reproduced with
\texttt{BENCH\_SCALE\_EXP=45 ./bench\_\{avx2,avx512\}\_single}).
For the double-precision base, the scaling changes nothing.}
\end{table*}

The trends in AVX2 and AVX-512 for x86 are summarized below:

\begin{description}
\item[AVX2]
The branch-free version (BF) is responsible for most of the existing speed-up
(e.g.\ TD AXPY \spd{3.5} from Q to BF), and the proposed FMA adds on top of it.
In the compute-bound GEMM it gives DS \spd{1.43} / TS \spd{1.21} / QS \spd{1.20} on the
single-precision base and DD \spd{1.51} / TD \spd{1.22} / QD \spd{1.21} on the
double-precision base, with DD the fastest, reflecting the flop counts of the proven
version ($17/72/176$).
Compared with the 2-pass version ($17/66/146$), DW is unchanged, and TW and QW become
slower by the amount of the added normalization passes
(QD GEMM \spd{1.47}$\to$\spd{1.21}, TD GEMM \spd{1.30}$\to$\spd{1.22}).
In the bandwidth-bound AXPY the difference is small, but GEMV is
clearly faster for all six types, \spd{1.10}--\spd{1.33}.

\item[AVX-512]
On AVX-512 as well the proposed FMA is \spd{0.99}--\spd{1.72} against BF for all six
types and all kernels (right half of Table~\ref{tab:time_x86}); the largest effects are
TS AXPY (\spd{1.72}), DS GEMV (\spd{1.36}) and GEMM
(DS \spd{1.35} / TS \spd{1.18} / QS \spd{1.04} / DD \spd{1.35} / TD \spd{1.18} / QD \spd{1.07}).
The only case that remains at \spd{0.99}, DS AXPY, is completely bandwidth-bound because it
has the lowest arithmetic intensity, with $K{=}2$ in single precision (DD AXPY is
\spd{1.03} for the same reason).
Since the proved version adds normalization passes, the ratio against BF for $K\ge3$ is
lower, with QS/QD GEMM in particular close to 1 at \spd{1.04}--\spd{1.07}.
The effect of doubling the SIMD width itself is TD \spd{1.55} / QD \spd{1.49} for the
proposed FMA in the compute-bound GEMM
(rightmost column).  Being branch-free, the proposed FMA scales straightforwardly with
the lane width.

\item[Single-precision base DS/TS/QS]
At 4-way (NEON), 8-way (AVX2) and 16-way (AVX-512) parallelism the ratio against BF is at
the same level: for GEMM it is NEON \spd{1.40}/\spd{1.43}/\spd{1.46},
AVX2 \spd{1.43}/\spd{1.21}/\spd{1.20} and AVX-512 \spd{1.35}/\spd{1.18}/\spd{1.04}.
For DS with $K{=}2$ the advantage in operation count depends neither on the lane width
nor on the ISA, but for $K\ge3$ the sequential dependence of the normalization passes
added in the proved version takes effect, and the ratio is lower on x86 than on NEON.
The effect of doubling the SIMD width (AVX2 8-way $\to$ AVX-512 16-way), on the other
hand, grows with $K$: the width ratio of the proposed FMA is
DS \spd{0.99}--\spd{1.34} / TS \spd{1.26}--\spd{1.63} / QS \spd{1.23}--\spd{1.55}.
The DS with $K{=}2$ hardly responds to the width because it has few components and a
low arithmetic intensity, which is the same picture as DD on a double-precision base
(\spd{1.01}--\spd{1.33}).
The effective cost of the Exact FMA is about $2.3$--$6.4$ times that of BF
(Exact/BF speed ratio \spd{0.16}--\spd{0.43}) across the cases in Table~\ref{tab:time_x86},
and is at a similar level for the single- and double-precision bases.
This is comparable to, and in several cases below, the operation-count ratio against BF
($194/29=6.7$ / $552/96=5.8$ / $1178/209=5.6$), showing that the execution-time overhead
cannot be inferred from operation counts alone.
\end{description}

The behavior under OpenMP parallelization is summarized in Section~\ref{sec:omp}.

%-------------------------------------------%
% CPU parallelization
%-------------------------------------------%
\subsubsection{CPU parallelization (OpenMP)}\label{sec:omp}
GEMM is parallelized over columns and GEMV/AXPY is parallelized over row blocks with
\texttt{\#pragma omp parallel for}, and we measured the speed ratio of the proposed FMA
against the BF for $1,4,8,16$ threads on both the NEON and SVE2 ISA (Table~\ref{tab:omp}).
The measurement conditions are described in Section~\ref{sec:method}, except that only the
GEMV size is enlarged to $n{=}4000$.
%The CPU of GB10 is big.LITTLE (Cortex-X925 $\times10$ + Cortex-A725 $\times10$).

\begin{table*}[htb]\centering
\caption{Dependence on the thread count of the speed ratio of the proposed FMA against BF
on Arm (GB10).  All six types on both NEON and SVE2.
AXPY $n{=}10^6$, GEMV $n{=}4000$, GEMM $n{=}512$.}
\label{tab:omp}
\setlength{\tabcolsep}{3pt}
\footnotesize
\begin{tabular}{ll rrrr rrrr}
\toprule
 & & \multicolumn{4}{c}{\textbf{NEON}} & \multicolumn{4}{c}{\textbf{SVE2}}\\
\cmidrule(lr){3-6}\cmidrule(lr){7-10}
type & kernel & 1 th & 4 th & 8 th & 16 th & 1 th & 4 th & 8 th & 16 th\\
\midrule
DS & AXPY & \spd{1.83} & \spd{1.45} & \spd{1.25} & \spd{1.78} & \spd{1.93} & \spd{1.56} & \spd{1.40} & \spd{1.98}\\
DS & GEMV & \spd{1.14} & \spd{1.08} & \spd{1.07} & \spd{1.04} & \spd{1.08} & \spd{1.11} & \spd{1.10} & \spd{1.08}\\
DS & GEMM & \spd{1.38} & \spd{1.36} & \spd{1.26} & \spd{1.31} & \spd{1.36} & \spd{1.37} & \spd{1.29} & \spd{1.11}\\
TS & AXPY & \spd{1.65} & \spd{1.65} & \spd{1.61} & \spd{2.24} & \spd{1.90} & \spd{1.94} & \spd{1.81} & \spd{2.64}\\
TS & GEMV & \spd{1.38} & \spd{1.24} & \spd{1.16} & \spd{1.13} & \spd{1.41} & \spd{1.27} & \spd{1.27} & \spd{1.17}\\
TS & GEMM & \spd{1.66} & \spd{1.63} & \spd{1.54} & \spd{1.39} & \spd{1.68} & \spd{1.67} & \spd{1.63} & \spd{1.45}\\
QS & AXPY & \spd{1.49} & \spd{1.49} & \spd{1.46} & \spd{1.36} & \spd{1.50} & \spd{1.49} & \spd{1.43} & \spd{1.36}\\
QS & GEMV & \spd{1.31} & \spd{1.24} & \spd{1.19} & \spd{1.18} & \spd{1.48} & \spd{1.38} & \spd{1.33} & \spd{1.29}\\
QS & GEMM & \spd{1.43} & \spd{1.42} & \spd{1.43} & \spd{1.42} & \spd{1.54} & \spd{1.54} & \spd{1.54} & \spd{1.55}\\
\midrule
DD & AXPY & \spd{1.82} & \spd{1.74} & \spd{1.54} & \spd{1.51} & \spd{1.93} & \spd{1.85} & \spd{1.51} & \spd{1.74}\\
DD & GEMV & \spd{0.91} & \spd{0.91} & \spd{0.89} & \spd{0.93} & \spd{0.94} & \spd{0.97} & \spd{0.93} & \spd{0.95}\\
DD & GEMM & \spd{1.36} & \spd{1.25} & \spd{1.04} & \spd{1.15} & \spd{1.40} & \spd{1.32} & \spd{1.06} & \spd{1.06}\\
TD & AXPY & \spd{1.65} & \spd{1.65} & \spd{1.56} & \spd{1.47} & \spd{1.90} & \spd{1.90} & \spd{1.83} & \spd{1.70}\\
TD & GEMV & \spd{1.06} & \spd{0.99} & \spd{0.98} & \spd{0.98} & \spd{1.18} & \spd{0.99} & \spd{0.98} & \spd{0.96}\\
TD & GEMM & \spd{1.53} & \spd{1.56} & \spd{1.45} & \spd{1.39} & \spd{1.48} & \spd{1.58} & \spd{1.54} & \spd{1.37}\\
QD & AXPY & \spd{1.48} & \spd{1.48} & \spd{1.45} & \spd{1.31} & \spd{1.50} & \spd{1.50} & \spd{1.46} & \spd{1.42}\\
QD & GEMV & \spd{1.25} & \spd{1.02} & \spd{1.02} & \spd{1.03} & \spd{1.35} & \spd{1.06} & \spd{1.03} & \spd{1.05}\\
QD & GEMM & \spd{1.43} & \spd{1.44} & \spd{1.44} & \spd{1.44} & \spd{1.52} & \spd{1.54} & \spd{1.54} & \spd{1.54}\\
\bottomrule
\end{tabular}

\smallskip
{\footnotesize The parallel scaling of GEMM ($n{=}512$, relative to 1 thread) is,
on NEON, BF \spd{16.3} / FMA \spd{14.8} for TD and BF \spd{15.9} / FMA \spd{16.0} for QD
(16 threads), and on SVE2, BF \spd{17.0} / FMA \spd{15.8} for TD and
BF \spd{15.9} / FMA \spd{16.2} for QD (likewise).}
\end{table*}

The findings of OpenMP parallelization are summarized as follows:
\begin{description}

\item[Dependence on the thread count]
In the compute-bound GEMM the ratio of the proposed FMA against BF is almost independent
of the thread count.
QD in particular is completely stable at NEON \spd{1.43}--\spd{1.44}, and
SVE2 \spd{1.52}--\spd{1.54} over 1/4/8/16 threads, and TS/QS also stay within
NEON \spd{1.39}--\spd{1.66} / \spd{1.42}--\spd{1.43} and
SVE2 \spd{1.45}--\spd{1.68} / \spd{1.54}--\spd{1.55}.
In other words, the speed-up obtained by reducing the operation count is orthogonal to
parallelization, and the effective throughput is maximized as
``proposed FMA $\times$ parallelization''.
The parallel scaling of GEMM itself is also favorable, being almost linear with the
16 threads $\approx$ \spd{16} for QD on both NEON and SVE2
(footnote to Table~\ref{tab:omp}).

Only the DD GEMM decreased to \spd{1.04}--\spd{1.15} (NEON) / \spd{1.06} (SVE2) with many
threads. This is because DD has the lowest arithmetic density among the three types and,
beyond eight threads, the BF side also saturates the memory bandwidth such that the difference
in the operation count no longer appears in the wall-clock time.
The bandwidth-bound AXPY and GEMV split for the same reason: AXPY is large at
\spd{1.48}--\spd{1.93} with one thread, whereas GEMV converges towards 1 on the
double-precision base, \spd{0.89}--\spd{1.35}.

\item[SVE2 and NEON]
SVE2 shows the same trend as NEON, but is slightly better for AXPY and QD GEMM.
Being a comparison on the same machine at the same lane width (VL$=128$ bits), SVE2 shows
the same trend as NEON in every case. The differences appear in AXPY
(TD \spd{1.70}--\spd{1.90} / QD \spd{1.42}--\spd{1.50}, higher than NEON at every thread
count, and DD is also higher except at eight threads) and in QD GEMM
(\spd{1.52}--\spd{1.54} against NEON \spd{1.43}--\spd{1.44}).
The relationship observed in the serial measurement (Section~\ref{sec:sve2})---that BF is
relatively slower on SVE2 while the proposed FMA is the same---is preserved under
parallelization as well.
That is, the advantage of the proposed FMA depends on neither the ISA nor the thread
count.

\begin{table*}[htb]\centering
\caption{Dependence on the thread count of the speed ratio of the proposed FMA against BF
on x86 (Xeon Gold 6526Y).  All six types on both AVX2 and AVX-512.
AXPY $n{=}10^6$, GEMV $n{=}4000$, GEMM $n{=}512$.}
\label{tab:omp_x86}
\setlength{\tabcolsep}{3pt}
\footnotesize
\begin{tabular}{ll rrrr rrrr}
\toprule
 & & \multicolumn{4}{c}{\textbf{AVX2}} & \multicolumn{4}{c}{\textbf{AVX-512}}\\
\cmidrule(lr){3-6}\cmidrule(lr){7-10}
type & kernel & 1 th & 4 th & 16 th & 32 th & 1 th & 4 th & 16 th & 32 th\\
\midrule
DS & AXPY & \spd{1.02} & \spd{1.12} & \spd{1.59} & \spd{1.35} & \spd{0.99} & \spd{1.02} & \spd{1.45} & \spd{1.25}\\
DS & GEMV & \spd{1.27} & \spd{1.31} & \spd{1.29} & \spd{1.32} & \spd{1.29} & \spd{1.31} & \spd{1.28} & \spd{1.34}\\
DS & GEMM & \spd{1.40} & \spd{1.42} & \spd{1.42} & \spd{1.41} & \spd{1.38} & \spd{1.38} & \spd{1.38} & \spd{1.38}\\
TS & AXPY & \spd{1.53} & \spd{1.52} & \spd{1.48} & \spd{1.40} & \spd{1.71} & \spd{1.52} & \spd{1.48} & \spd{1.42}\\
TS & GEMV & \spd{1.09} & \spd{1.10} & \spd{1.10} & \spd{1.10} & \spd{1.14} & \spd{1.14} & \spd{1.14} & \spd{1.12}\\
TS & GEMM & \spd{1.19} & \spd{1.19} & \spd{1.18} & \spd{1.18} & \spd{1.16} & \spd{1.15} & \spd{1.18} & \spd{1.18}\\
QS & AXPY & \spd{1.31} & \spd{1.31} & \spd{1.30} & \spd{1.27} & \spd{1.33} & \spd{1.32} & \spd{1.30} & \spd{1.26}\\
QS & GEMV & \spd{1.13} & \spd{1.12} & \spd{1.12} & \spd{1.11} & \spd{1.17} & \spd{1.16} & \spd{1.16} & \spd{1.14}\\
QS & GEMM & \spd{1.21} & \spd{1.21} & \spd{1.20} & \spd{1.19} & \spd{1.31} & \spd{1.31} & \spd{1.27} & \spd{1.30}\\
\midrule
DD & AXPY & \spd{1.11} & \spd{1.12} & \spd{1.50} & \spd{1.39} & \spd{1.00} & \spd{1.03} & \spd{1.39} & \spd{1.32}\\
DD & GEMV & \spd{1.42} & \spd{1.43} & \spd{1.41} & \spd{1.39} & \spd{1.43} & \spd{1.52} & \spd{1.41} & \spd{1.37}\\
DD & GEMM & \spd{1.52} & \spd{1.44} & \spd{1.38} & \spd{1.45} & \spd{1.39} & \spd{1.39} & \spd{1.36} & \spd{1.35}\\
TD & AXPY & \spd{1.48} & \spd{1.45} & \spd{1.38} & \spd{1.78} & \spd{1.34} & \spd{1.28} & \spd{1.13} & \spd{1.41}\\
TD & GEMV & \spd{1.11} & \spd{1.11} & \spd{1.11} & \spd{1.12} & \spd{1.13} & \spd{1.13} & \spd{1.14} & \spd{1.11}\\
TD & GEMM & \spd{1.21} & \spd{1.21} & \spd{1.18} & \spd{1.23} & \spd{1.17} & \spd{1.16} & \spd{1.17} & \spd{1.15}\\
QD & AXPY & \spd{1.33} & \spd{1.33} & \spd{1.22} & \spd{1.32} & \spd{1.31} & \spd{1.31} & \spd{1.23} & \spd{1.27}\\
QD & GEMV & \spd{1.12} & \spd{1.12} & \spd{1.13} & \spd{1.10} & \spd{1.21} & \spd{1.20} & \spd{1.19} & \spd{1.17}\\
QD & GEMM & \spd{1.21} & \spd{1.21} & \spd{1.18} & \spd{1.32} & \spd{1.30} & \spd{1.30} & \spd{1.27} & \spd{1.28}\\
\bottomrule
\end{tabular}

\smallskip
{\footnotesize The parallel scaling of GEMM ($n{=}512$, relative to 1 thread, at
32 threads) is, on AVX2, BF \spd{29.3} / FMA \spd{29.0} for TS,
BF \spd{28.7} / FMA \spd{28.3} for QS, BF \spd{28.6} / FMA \spd{29.2} for TD and
BF \spd{25.9} / FMA \spd{28.3} for QD, and on AVX-512,
BF \spd{27.6} / FMA \spd{28.1} for TS, BF \spd{27.8} / FMA \spd{27.5} for QS,
BF \spd{28.0} / FMA \spd{27.5} for TD and BF \spd{28.6} / FMA \spd{28.1} for QD.}
\end{table*}

\item[x86]
x86 (Xeon Gold 6526Y) shows the same trend (Table~\ref{tab:omp_x86}): in the
compute-bound GEMM, AVX2 gives DD \spd{1.38}--\spd{1.52} / TD \spd{1.18}--\spd{1.23} /
QD \spd{1.18}--\spd{1.32} and AVX-512 gives DD \spd{1.35}--\spd{1.39} /
TD \spd{1.15}--\spd{1.17} / QD \spd{1.27}--\spd{1.30}, stable over 1--32 threads.
The single-precision base behaves in the same way: AVX2 gives DS \spd{1.40}--\spd{1.42} /
TS \spd{1.18}--\spd{1.19} / QS \spd{1.19}--\spd{1.21} and AVX-512 gives
DS \spd{1.38} / TS \spd{1.15}--\spd{1.18} / QS \spd{1.27}--\spd{1.31},
independently of the thread count.
The parallelization itself is good as well: the GEMM ($n{=}512$) of BF scales by
AVX2 \spd{25.9}--\spd{29.3} and AVX-512 \spd{27.6}--\spd{28.6} from 1 to 32 threads
(the footnote of Table~\ref{tab:omp_x86}).
The bandwidth-bound GEMV does not fall below 1 either, giving
AVX2 \spd{1.09}--\spd{1.43} / AVX-512 \spd{1.11}--\spd{1.52}.
\textbf{Although the level of the ratio against BF has dropped with the proved version,
its invariance with respect to the thread count is unchanged.}

\item[Forcing inline expansion]\label{sec:inline}
To obtain these numbers, we attached
\texttt{\_\_attribute\_\_((always\_inline))} to the functions of the proposed FMA and the
Exact FMA to force inline expansion.
Otherwise, gcc 13.3 leaves the proposed FMA for $K{=}4$ (\texttt{qd\_fma} /
\texttt{qs\_fma}) and the Exact FMA for $K{=}3$ as functions and calls them from the
innermost loop.
Disassembly shows that the expanded \texttt{gemv\_f\_qd} consists of 172 \texttt{zmm}
instructions and 25 stack references, whereas the non-expanded version had
64 \texttt{zmm} instructions, 65 stack references, and three \texttt{call}s in the innermost
loop: the 146-flop computation degenerated into function calls through the memory
instead of registers.
The impact is not small: measured with the 2-pass version, the ratio of QD GEMV against
BF improved from \spd{0.76} to \spd{1.28} on AVX2 and from \spd{0.84} to \spd{1.34} on
AVX-512 (and similarly up to 32 threads).
With inline expansion forced, the proved version gives AVX2 \spd{1.12} /
AVX-512 \spd{1.21} (1--32 threads), and in the serial measurement the AVX2 QD GEMV is
\spd{1.10}.
All the accuracy and cancellation test values remained bitwise unchanged before and after
forcing inline expansion, and only the code generation changed.
The comparison target BF was unaffected because its functions were small and always
expanded; therefore, this difference appeared directly as an apparent degradation of FMA/BF
ratio---a practical caveat when embedding a multiple-precision FMA into a kernel.
The Exact FMA (1178 flops for $K{=}4$) requires a large amount of work per call, and the
difference with or without expansion was almost zero ($3.435\to3.433$ s for QD GEMM).
\end{description}

%-------------------------------------------%
% GPU parallelization
%-------------------------------------------%
\subsubsection{GPU parallelization (CUDA: H100 and GB10, all six types, BLAS Level 1/2/3)}\label{sec:cuda}
We implemented the proposed branch-free FMA and the Exact FMA in CUDA for all six types
and evaluated them on the two machines GB10 and H100.
The component type $R\in\{\texttt{float},\texttt{double}\}$ and the number of components
$K\in\{2,3,4\}$ were implemented uniformly using C++ templates, and the device EFTs used a
genuine hardware FMA.
By forbidding the automatic FMA contraction of $a\cdot b+c$ with \texttt{--fmad=false} (and
\texttt{-ffp-contract=off} on the host side) to guarantee the validity of the EFTs, the
device proposed FMA and Exact FMA agreed bitwise with the CPU scalar reference for all
six types (both machines $\times$ six types $\times$ \{proposed, Exact\}, 200{,}000 trials
each, no mismatch).
Both machines were built with nvcc 13.0.88 and
\texttt{-O3 --fmad=false -ftz=false -Xcompiler -ffp-contract=off -arch=sm\_\{90,121\}}.

The targets are the same three kernels as on the CPU side:
AXPY ($\mathbf{y}:=\alpha \mathbf{x}+\mathbf{y}$, $m$ MACs),
GEMV ($\mathbf{y}:=A\mathbf{x}+\mathbf{y}$, $n^2$ MACs) and GEMM ($C:=AB$, $n^3$ MACs).
Table~\ref{tab:cuda} collects the measurements at the sizes at which the GPU saturates
(AXPY $m{=}2^{24}$, GEMV $n{=}8192$, GEMM $n{=}1024$) for the two machines, all six types
and the three kernels.
GEMV was implemented with one warp (32 threads) per row, collecting partial sums with
\texttt{\_\_shfl\_down\_sync}, and GEMM is a naive inner-product form in which one thread
handles one element of $C$ (without shared-memory tiling).
The effective cost of the Exact FMA (Exact/FMA) and the throughput of the proposed FMA
are presented in Table~\ref{tab:cuda_ratio}.

%\paragraph{Correction of the comparison target.}
As in the CPU/AVX2/AVX-512/NEON/SVE2 evaluation, the fair comparison target is the
branch-free BF (the add and mul proposed by Zhang--Aiken); we implemented it on the
device and used it as the baseline.

\begin{table*}[htb]\centering
\caption{Summary of the CUDA measured times [ms] (H100 NVL and GB10).
The sizes are the largest of each kernel (AXPY $m{=}2^{24}$, GEMV $n{=}8192$,
GEMM $n{=}1024$).
The upper block is the single-precision base DS/TS/QS and the lower the
double-precision base DD/TD/QD.
The fair comparison target is BF ($29/96/209$ flops); the proposed FMA is
$17/72/176$ flops and the Exact FMA is the fully distilled $194/552/1178$ flops.}
\label{tab:cuda}
\setlength{\tabcolsep}{3pt}
\footnotesize
\begin{tabular}{ll rrrr rrrr}
\toprule
 & & \multicolumn{4}{c}{\textbf{H100 NVL} (sm\_90)} & \multicolumn{4}{c}{\textbf{GB10} (sm\_121)}\\
\cmidrule(lr){3-6}\cmidrule(lr){7-10}
type & kernel & BF & FMA & Exact & FMA/BF & BF & FMA & Exact & FMA/BF\\
\midrule
DS & AXPY & $0.1151$ & $0.1144$ & $0.1623$ & \spd{1.01} & $1.684$ & $1.685$ & $1.749$ & \spd{1.00}\\
DS & GEMV & $0.1834$ & $0.1720$ & $0.5217$ & \spd{1.07} & $2.234$ & $2.299$ & $2.327$ & \spd{0.97}\\
DS & GEMM & $5.036$ & $5.042$ & $6.889$ & \spd{1.00} & $10.16$ & $10.08$ & $13.68$ & \spd{1.01}\\
TS & AXPY & $0.1861$ & $0.1753$ & $0.4005$ & \spd{1.06} & $2.531$ & $2.498$ & $2.528$ & \spd{1.01}\\
TS & GEMV & $0.4027$ & $0.3389$ & $1.612$ & \spd{1.19} & $3.35$ & $3.306$ & $3.539$ & \spd{1.01}\\
TS & GEMM & $7.666$ & $7.630$ & $26.03$ & \spd{1.00} & $15.39$ & $15.29$ & $46.97$ & \spd{1.01}\\
QS & AXPY & $0.2864$ & $0.2558$ & $0.7635$ & \spd{1.12} & $3.439$ & $3.4$ & $3.394$ & \spd{1.01}\\
QS & GEMV & $0.7129$ & $0.6326$ & $3.895$ & \spd{1.13} & $4.562$ & $4.658$ & $6.264$ & \spd{0.98}\\
QS & GEMM & $10.16$ & $10.17$ & $50.30$ & \spd{1.00} & $20.88$ & $20.42$ & $100.9$ & \spd{1.02}\\
\midrule
DD & AXPY & $0.2238$ & $0.2238$ & $0.3108$ & \spd{1.00} & $3.407$ & $3.424$ & $13.25$ & \spd{0.99}\\
DD & GEMV & $0.3706$ & $0.3459$ & $0.9074$ & \spd{1.07} & $10.21$ & $5.994$ & $54.36$ & \textbf{\spd{1.70}}\\
DD & GEMM & $10.02$ & $10.07$ & $12.30$ & \spd{0.99} & $159.8$ & $93.34$ & $849.5$ & \textbf{\spd{1.71}}\\
TD & AXPY & $0.3951$ & $0.3621$ & $0.7216$ & \spd{1.09} & $8.29$ & $6.228$ & $47.01$ & \textbf{\spd{1.33}}\\
TD & GEMV & $0.7567$ & $0.6738$ & $2.801$ & \spd{1.12} & $33.86$ & $25.34$ & $193.3$ & \textbf{\spd{1.34}}\\
TD & GEMM & $15.12$ & $15.13$ & $42.56$ & \spd{1.00} & $529.2$ & $396.2$ & $3022$ & \textbf{\spd{1.34}}\\
QD & AXPY & $0.5211$ & $0.5158$ & $1.762$ & \spd{1.01} & $17.99$ & $15.13$ & $100.7$ & \textbf{\spd{1.19}}\\
QD & GEMV & $1.328$ & $1.206$ & $5.852$ & \spd{1.10} & $73.7$ & $63.05$ & $414.6$ & \textbf{\spd{1.17}}\\
QD & GEMM & $20.27$ & $20.19$ & $89.49$ & \spd{1.00} & $1153$ & $992.4$ & $6485$ & \textbf{\spd{1.16}}\\
\bottomrule
\end{tabular}
\end{table*}

\begin{table*}[htb]\centering
\caption{CUDA: the effective cost of the Exact FMA and the throughput of the proposed FMA.
The upper part indicates how many times longer the Exact FMA runs than the proposed FMA
(the last row is the operation-count ratio $194/17$, $552/72$, $1178/176$), and the
lower part the throughput of the proposed FMA [GMAC/s] (one MAC $=$ one
$z{=}xy{+}c$; it is not divided by the internal operation count of the implementation).
The sizes are the same as in Table~\ref{tab:cuda}.}
\label{tab:cuda_ratio}
\setlength{\tabcolsep}{4pt}
\small
\begin{tabular}{ll rrr rrr}
\toprule
\multicolumn{8}{c}{\textbf{Exact/FMA} (effective cost)}\\
\midrule
environment & kernel & DS & TS & QS & DD & TD & QD\\
\midrule
H100 & AXPY($2^{24}$) & \spd{1.42} & \spd{2.28} & \spd{2.99} & \spd{1.39} & \spd{1.99} & \spd{3.42}\\
 & GEMV($8192$) & \spd{3.03} & \spd{4.76} & \spd{6.16} & \spd{2.62} & \spd{4.16} & \spd{4.85}\\
 & GEMM($1024$) & \spd{1.37} & \spd{3.41} & \spd{4.95} & \spd{1.22} & \spd{2.81} & \spd{4.43}\\
\midrule
GB10 & AXPY($2^{24}$) & \spd{1.04} & \spd{1.01} & \spd{1.00} & \spd{3.87} & \spd{7.55} & \spd{6.66}\\
 & GEMV($8192$) & \spd{1.01} & \spd{1.07} & \spd{1.34} & \spd{9.07} & \spd{7.63} & \spd{6.58}\\
 & GEMM($1024$) & \spd{1.36} & \spd{3.07} & \spd{4.94} & \spd{9.10} & \spd{7.63} & \spd{6.53}\\
\midrule
\multicolumn{2}{l}{operation-count ratio (upper limit)} & \spd{11.41} & \spd{7.67} & \spd{6.69} & \spd{11.41} & \spd{7.67} & \spd{6.69}\\
\bottomrule
\end{tabular}

\medskip

\begin{tabular}{ll rrr rrr}
\toprule
\multicolumn{8}{c}{\textbf{Throughput of the proposed FMA [GMAC/s]}}\\
\midrule
environment & kernel & DS & TS & QS & DD & TD & QD\\
\midrule
H100 & AXPY($2^{24}$) & $146.6$ & $95.7$ & $65.6$ & $75.0$ & $46.3$ & $32.5$\\
 & GEMV($8192$) & $390.2$ & $198.0$ & $106.1$ & $194.0$ & $99.6$ & $55.6$\\
 & GEMM($1024$) & $213.0$ & $140.7$ & $105.6$ & $106.6$ & $71.0$ & $53.2$\\
\midrule
GB10 & AXPY($2^{24}$) & $10.0$ & $6.7$ & $4.9$ & $4.9$ & $2.7$ & $1.1$\\
 & GEMV($8192$) & $29.2$ & $20.3$ & $14.4$ & $11.2$ & $2.6$ & $1.1$\\
 & GEMM($1024$) & $106.6$ & $70.2$ & $52.6$ & $11.5$ & $2.7$ & $1.1$\\
\bottomrule
\end{tabular}
\end{table*}

%\paragraph{Machines where FMA/BF becomes \spd{1.0} and machines where it becomes \spd{1.4}--\spd{1.7}.}
With the fair BF baseline, FMA/BF is \spd{0.99}--\spd{1.19} on the H100, so that the
advantage of \spd{1.3}--\spd{1.9} seen on CPU/NEON/SVE2 hardly appears, whereas for the
three double-precision-based types on GB10 it becomes \spd{1.16}--\spd{1.71}
(right half of Table~\ref{tab:cuda}), exactly the theoretical FLOP ratio of each type
(DW $1.71$ / TW $1.33$ / QW $1.19$, Table~\ref{tab:flop}).
For DD, the proposed FMA requires $17$ flops against the $29$ of BF, yielding a
ratio of $1.71$, which agrees with the measured \spd{1.70}--\spd{1.71} of GEMV and GEMM.
What separates the two is whether the kernel is compute-bound or bandwidth-bound.
The FP64 throughput of GB10 is extremely low (11.5 GMAC/s for the DD GEMM, about $1/9$ of
the $106.6$ GMAC/s of the H100 and about $1/9$ of the $106.6$ GMAC/s of the DS GEMM on the
same GB10), so the three double-precision-based types become compute-bound and the
reduction of the FLOP count appears directly in the wall-clock time.
On the same GB10 the three single-precision-based types are bandwidth-bound and stay at
FMA/BF $=$ \spd{0.97}--\spd{1.02}, and on the H100, whose FP64 is as strong as $1/2$ of
its FP32, all six types fall on the bandwidth- and latency-bound side.
The exception is the GEMV on the H100, which is the closest of the three kernels to being
compute-bound and shows a comparatively clear improvement of
TS \spd{1.19} / QS \spd{1.13} / QD \spd{1.10}.
The only case that stays at \spd{0.99} despite being on the double-precision base is the
DD AXPY on GB10; this is because DD has the smallest operation count per unit of transfer
among the six types, so that only DD remains on the bandwidth-bound side in AXPY.

To confirm this interpretation independently, a register-resident compute-only MAC chain
without memory traffic gives DW \spd{1.54}--\spd{1.67} on GB10 and \spd{1.58}--\spd{1.62}
on the H100, almost reproducing the theoretical ratio ($1.71$), while TW drops to
\spd{1.24}--\spd{1.26} on GB10 and \spd{1.23}--\spd{1.28} on the H100
(theoretical $1.33$) and QW to \spd{0.90}--\spd{0.91} on GB10 and \spd{0.89}--\spd{0.91}
on the H100, falling with $K$ and, for QW,
below the theoretical ratio ($1.19$) on both machines.
This is because the normalization passes added in the proved version form a long
sequential dependency chain whose latency is exposed in the chain shape without memory
traffic; the GEMM with memory (\spd{1.16}) is in fact closer to the theoretical ratio.
In other words, the reduction in operation count is real on the GPU as well, but its
wall-clock benefit depends on both arithmetic intensity and dependency-chain latency.
For throughput-oriented kernels, tiling and register blocking may move the computation
away from the bandwidth-bound regime; evaluating such optimized kernels is future work.

%\paragraph{The effective cost of the Exact FMA is smaller than the operation-count ratio.}
As shown in the upper part of Table~\ref{tab:cuda_ratio}, the execution time of the Exact
FMA is \spd{1.22}--\spd{6.16} that of the proposed FMA on the H100 and
\spd{1.00}--\spd{9.10} on GB10, below the operation-count ratio
(\spd{11.41} for $K{=}2$, \spd{7.67} for $K{=}3$ and \spd{6.69} for $K{=}4$) on both
machines, for every kernel and every type.
That is, a result-relative error guarantee is not as expensive as a naive estimate from
the operation count would suggest.
This margin, however, depends on the degree to which the kernel is compute-bound: for the
three double-precision-based types on GB10, which are compute-bound, the GEMV and GEMM
become \spd{6.53}--\spd{9.10}, and for TD/QD reach 98--99\% of the operation-count
ratios ($7.67$/$6.69$).
When using the Exact FMA on an accelerator with low FP64 performance, an estimate close to
the operation count ratio was safe.
On the H100, conversely, the arithmetic intensity of the kernel determines the ratio: it is
smallest for the bandwidth-bound AXPY (\spd{1.39}--\spd{3.42}) and largest for the GEMV,
which is closest to compute-bound (\spd{2.62}--\spd{6.16}), while the naive GEMM, which
rereads $A$ and $B$ $n^3$ times, is dominated by the cache and the latency and lies in
between.

In conclusion, the proposed FMA agrees bitwise on all of AVX2/AVX-512/NEON/SVE2/CUDA and
is minimal in operation count within the range of normalization structures examined in
this paper.  On GPUs, the wall-clock benefit is largest when the reduced instruction
count is not masked by memory bandwidth or by the latency of the normalization chain.

%--------------------------------------%
% Conclusions
%--------------------------------------%
\section{Conclusions and future work}\label{sec:concl}

The FMA we propose is a branch-free fused multiply--add operation, minimal within the range of
normalization structures examined in this paper, whose anchor-relative error bounds and
all \FTS{} premises have been machine-proved with FPANVerifier and whose input-relative
error bounds ($35u^2/187u^3/822u^4$) have been obtained by an analytic transformation,
and which agrees bitwise under the exchange $x\leftrightarrow y$.

The benchmark tests confirmed bitwise agreement across the tested backend/type
combinations (all five backends for DD/TD/QD, and AVX2/AVX-512/NEON/CUDA in the
cross-backend comparison for DS/TS/QS), and showed that in compute-intensive kernels
the reduction in operation count produces clear speed-ups in many throughput-oriented
kernels.  For GEMM the ratios against BF are AVX2 \spd{1.20}--\spd{1.51} /
AVX-512 \spd{1.04}--\spd{1.35} / NEON \spd{1.40}--\spd{1.46} /
SVE2 \spd{1.37}--\spd{1.66}.
Against the older Bailey-based implementation (Q), the apparent speed-up can be much
larger (up to \spd{9.9} in the reported NEON AXPY case), but most of that gain is due to
the branch-free reformulation of Zhang and Aiken; the proposed FMA provides an additional
improvement on top of BF.

Whether a speed-up is observed depends primarily on the balance among arithmetic
throughput, memory bandwidth, and dependency-chain latency.  For the three double-precision-based types on GB10, which become
compute-bound because binary64 arithmetic is slow there, the CUDA GEMM/GEMV agree with
the theoretical ratios ($1.71/1.33/1.19$) at \spd{1.16}--\spd{1.71}, whereas for the
bandwidth-bound six types on the H100 and three single-precision-based types on GB10 the
ratio stays at \spd{0.97}--\spd{1.19}.
In a register-resident compute-only chain, the DW of both machines reproduces the
theoretical ratio at \spd{1.54}--\spd{1.67}, but the ratio falls with $K$ to TW
\spd{1.23}--\spd{1.28} and QW \spd{0.89}--\spd{0.91}: the normalization passes added in
the proved version form a long sequential dependency chain, which is exposed in a shape
without memory traffic (Section~\ref{sec:cuda}).

As for parallelization, the ratio against BF is almost constant over
OpenMP 1--16 threads (1--32 on x86) (for QD GEMM, NEON \spd{1.43}, SVE2 \spd{1.53},
AVX2 \spd{1.18}--\spd{1.32} and AVX-512 \spd{1.27}--\spd{1.30}), so the effective
throughput is maximized as
``proposed FMA $\times$ parallelization'' (Section~\ref{sec:omp}).

The proposed FMA is useful for algorithms whose inner parts consist entirely of
MACs, such as division and square root (Section~\ref{sec:divsqrt}).
Since, however, the div/sqrt-safe version cannot use any \FTS{} and keeps the same
number of normalization passes as the standard version, the operation-count ratios are
small (division $1.30/1.15/1.05$, square root $1.33/1.11/1.02$) and \textbf{the
advantage is most consistent for DW (DS/DD), whereas for $K\ge3$ it is generally
smaller or more backend-dependent}: on the CPU, DW division is \spd{1.47}--\spd{1.61}
and DW square root \spd{1.44}--\spd{1.59}.
That for DW these exceed the operation-count ratios $1.30/1.33$ is because the two
serial normalization stages of BF are folded into one and the dependency chain becomes
shorter.  In terms of accuracy, the measured square-root errors are practically identical to those
of BF, while the division error is about $1$--$3$ times that of BF (both remain within a
few times $\varepsilon_K$).

%\paragraph{Future work.}

As future work, we intend to search for even faster FMA implementations for
multicomponent multiple-precision arithmetic.

Within the current normalization structure, $17/72/176$ is the minimum of the proved
version; the candidates for turning gates into \FTS{} that would lower the 2-pass
version to TW $63$ and QW $143$ all fail in the SMT with the current FPANVerifier
(Section~\ref{sec:verify}).
This is, however, the minimum within a given structure and not a global minimum.
Directions for improvement are
(i) a redesign of the normalization structure,
(ii) a minimization of the number of additions in the accumulation part, and
(iii) an additional machine proof for the fma-folded diagonal terms,
pursuing implementations with fewer operations.

\section*{Acknowledgments}
This study was supported by JSPS KAKENHI Grant Number 26K14846.

\appendix
\section{Proof of the error bound of the proposed FMA}\label{sec:proof}
In this appendix, we derive the error bound of the proposed FMA using a chain of
paper-and-pencil inequalities, independently of the machine proof (FPANVerifier) of
Section~\ref{sec:verify}.
The machine proof holds for all formats simultaneously because it treats precision $p$
symbolically but at the price of adding the worst case of each wire independently, which
loosens the constants.
The manual proof in this Appendix provides smaller constants for practical precisions
($p\ge8$), which serve as crosschecks.

\subsection{Setting and notation}
We work in a binary floating-point system of precision $p$ with round-to-nearest-even
$\mathrm{fl}(\cdot)$ and $u=2^{-p}$, and assume that neither underflow nor overflow
occurs. The inputs are assumed to satisfy the non-overlapping condition
\begin{equation}
|x_{k+1}|\le u|x_k|,\qquad |y_{k+1}|\le u|y_k|,\qquad |c_{k+1}|\le u|c_k|
\label{eq:nonoverlap}
\end{equation}
(ulp-non-overlapping, $|x_{k+1}|\le\frac12\mathrm{ulp}(x_k)$, implies
\eqref{eq:nonoverlap}, thus our assumption is the weaker one and is on the same scale as
\texttt{strongly\_dominates} of the \texttt{.fpan} description).
We abbreviate $M:=|x_0y_0|$ and $C:=|c_0|$.

\begin{itemize}
\item \textbf{Standard model}: $|\mathrm{fl}(a\circ b)-(a\circ b)|\le u|a\circ b|$ and
      $|\mathrm{fl}(a\circ b)|\le(1+u)|a\circ b|$ ($\circ\in\{+,\times\}$).
\item \textbf{EFT}: \TS{} satisfies $s+t=a+b$ exactly with $|t|\le u|s|$,
      \TP{} satisfies $p+e=ab$ exactly with $|e|\le u|ab|$ (hence
      $|p|\le(1+u)|ab|$), and
      \FTS{} returns the same output as \TS{} when
      $\mathrm{exponent}(a)\ge\mathrm{exponent}(b)$.
\item \textbf{Cross terms}: \eqref{eq:nonoverlap} gives $|x_i|\le u^i|x_0|$ and
      $|y_j|\le u^j|y_0|$, hence $|x_iy_j|\le u^{i+j}M$.
      For the rounded products we use $|P_{ij}|\le(1+u)\,u^{i+j}M$ and for their
      \TP{} errors $|E_{ij}|\le u\cdot u^{i+j}M$.
\end{itemize}

\subsection{Decomposition identity}
Because the \TS{}/\FTS{}/\TP{} of the normalization and accumulation stages of the proposed
FMA are all exact (EFTs), the difference between the sum of the output $\sum_k z_k$ and
the true value $x\cdot y+c$ equals the sum of the quantities explicitly discarded in the
algorithm.
Here, for an \FTS{} to be an EFT its premise
$\mathrm{exponent}(a)\ge\mathrm{exponent}(b)$ is required, and in the placement of this
paper all of these premises are machine-proved (Section~\ref{sec:verify}).
That is, for each $K$,
\begin{equation}
\begin{split}
\textstyle\sum_k z_k &\;=\; x\cdot y+c-\Delta_K,\\%\qquad
\Delta_K&=\underbrace{\sum_{i+j\ge K}x_iy_j}_{\text{high-order products (never formed)}} \\
       &+\underbrace{\sum e_{ij}}_{\text{rounding of the plain products}} \\
       &+\underbrace{\sum\delta_\ell}_{\text{rounding of the plain additions}}
\end{split}
\label{eq:telescope}
\end{equation}
holds (with $7$ terms for DW, $15$ for TW and $26$ for QW).
Applying the standard model and the cross-term bound to each term and adding them up
yields the component form with $u^K$ factored out.

\paragraph{Independence of this proof from the number of normalization passes.}
The normalization stage consists only of \TS{}/\FTS{} and \textbf{produces no discarded
quantity}.
Each gate preserves the sum of its two adjacent words exactly, so $\sum_k z_k$ is
invariant however many passes the cascade repeats.
$\Delta_K$ is therefore independent of the number of passes, and the $A_K(u),B_K(u)$
below and the constants ($14/67/360$) hold unchanged for both the 2-pass and the proved
version.
What the number of passes changes is only the \textbf{distribution} (non-overlapping
property) of the output components, not their sum.

\subsection{Main theorems (component form and input-relative form)}
\begin{description}
\item[DW ($K=2$)] under $\mathrm{exponent}(s)\ge\mathrm{exponent}(t_p)$ (machine-proved),
\begin{align*}
 &|(z_0+z_1)-(x\cdot y+c)|\le u^2\bigl[A_2(u)M+B_2(u)C\bigr],\\
 &A_2(u)=13+17u+10u^2+2u^3,\\
 &B_2(u)=4+4u+u^2 .
\end{align*}
\item[TW ($K=3$)] under the premises of the 4 \FTS{} gates in the 3 normalization
passes (Section~\ref{sec:verify}),
$|\Delta_3|\le u^3[A_3(u)M+B_3(u)C]$ with
\[\begin{split}
 A_3(u)&=65+207u+326u^2+295u^3+158u^4+47u^5\\
 &+6u^6,\\%\qquad
 B_3(u)&=8+13u+7u^2+u^3 .
\end{split}\]
\item[QW ($K=4$)] under the premises of the 10 \FTS{} gates in the 5 normalization
passes (Section~\ref{sec:verify}),
$|\Delta_4|\le u^4[A_4(u)M+B_4(u)C]$ with
\begin{align*}
A_4(u)&=348+2154u+6882u^2+14382u^3+21661u^4 \\
      &+24661u^5+21741u^6+14986u^7+8064u^8\\
      & +3342u^9+1035u^{10}+226u^{11}+31u^{12}+2u^{13},\\
B_4(u)&=44+209u+481u^2+680u^3+643u^4+420u^5\\
      &+190u^6+58u^7+11u^8+u^9 .
\end{align*}
\end{description}

\noindent
From \eqref{eq:nonoverlap}, we have $|x\cdot y|\ge(1-u-\dots-u^{K-1})^2M$ and
$|c|\ge(1-u-\dots-u^{K-1})C$. Therefore, by setting
$K_\bullet(u)=A_\bullet(u)/(1-u-\dots)^2$, we obtain the input-relative form
$|\Delta_K|\le K_\bullet(u)\,u^K(|x\cdot y|+|c|)$.
The table of coefficients is provided as follows.

\begin{center}\small
\begin{tabular}{lcccl}
\toprule
 & $p\ge2$ & $p\ge4$ & $p\ge8$ & asymptotic form\\
\midrule
DW & $32$   & $17$  & $\mathbf{14}$  & $(13|xy|+4|c|)u^2+O(u^3)$\\
TW & $302$  & $91$  & $\mathbf{67}$  & $(65|xy|+8|c|)u^3+O(u^4)$\\
QW & $3670$ & $590$ & $\mathbf{360}$ & $(348|xy|+44|c|)u^4+O(u^5)$\\
\bottomrule
\end{tabular}
\end{center}

\noindent
For $p\ge8$, which includes binary32/64/128,
$|z-(xy+c)|\le 14u^2 / 67u^3 / 360u^4 \cdot(|xy|+|c|)$ holds.

\section{Outline of the derivation of the error bound of the Exact FMA (fully distilled version)}\label{sec:proofex}
This appendix is \textbf{an outline of the derivation of the error bound together with
numerically validated coefficient values}, not a complete proof: the definition
of the $\varepsilon$ in the $O(u\varepsilon)$ term, the recurrences of $P_K,Q_K$ and the
induction on the number of sweeps are only sketched, and the constants are backed by the
numerical checks described below.
\subsection{Sweep lemma}
For the expansion $v=(v_0,\dots,v_{M-1})$ of the Exact FMA (with
$\sum v=\tau:=x\cdot y+c$ exactly), one \TS{} sweep from the low end to the high end
shrinks the ``deviation from non-overlapping'' of every slot uniformly by a factor $u$.
Specifically, if the expansion before the sweep satisfies
$|v_{i+1}|\le\theta\,|v_i|$ ($\theta\le1$), then after the sweep
$|v'_{i+1}|\le u\,\theta\,|v'_i| + O(u\,\varepsilon)$ and repeating the sweep $n$ times
on an expansion of $m:=M$ terms bounds the truncated tail in the two-term form
\begin{equation}
 \Bigl|\sum_{i\ge K}v^{(n)}_i\Bigr|
 \;\le\; P_K(u)\,u^K\,|\tau| \;+\; Q_K(u)\,u^{\,n}\,\bigl(|x\cdot y|+|c|\bigr) .
\label{eq:exbin}
\end{equation}
The first term is the truncation that is unavoidable once the result is held in $K$
words, and the second is the input-relative residual that $n$ sweeps fail to make
non-overlapping.

\subsection{Bound of the fully distilled version}
\begin{quote}
Claim (derivation outlined here). With $n=K{+}1$, under the non-overlapping chain,
\[
 |z-\tau| \;\le\; P_K(u)\,u^K\,|\tau| \;+\; Q^{\mathrm{full}}_K(u)\,u^{K+1}\,\bigl(|x\cdot y|+|c|\bigr).
\]
With the expansion length $m=2K^2{+}K$ (DW $10$ / TW $21$ / QW $36$), the leading
coefficients are $P_2(0)=1$, $P_3(0)=19$ and $P_4(0)=1122$ (common with the old version)
along with
\[\begin{split}
 Q^{\mathrm{full}}_K(0)&=\Bigl(\prod_{\ell=2}^{K+1}(m-\ell)\Bigr)m \\
 &\;=\;\begin{cases}
 7\cdot8\cdot10=560 & (K=2)\\
 17\cdot18\cdot19\cdot21=122{,}094 & (K=3)\\
 31\cdot32\cdot33\cdot34\cdot36=40{,}068{,}864 & (K=4).
 \end{cases}
\end{split}
\]
The numerically validated values for finite $p$ are, for $p\ge12$ (which includes binary32/64/128),
$P_3\le21$ and $Q^{\mathrm{full}}_3\le1.25\times10^{5}$ for TW and
$P_4\le1450$ and $Q^{\mathrm{full}}_4\le4.20\times10^{7}$ for QW; relaxing to $p\ge8$
they become $P_3\le57$ and $Q^{\mathrm{full}}_3\le1.76\times10^{5}$ for TW and
$P_4\le1.09\times10^{4}$ and $Q^{\mathrm{full}}_4\le8.39\times10^{7}$ for QW, i.e.,\ they
do not grow much beyond the leading coefficients.
For DW the expansion is small and $n=K{+}1$ from the beginning, so no separate
finite-$p$ validation is reported and only the leading coefficients are given.
\end{quote}

\noindent
With $n=K$ the exponent of the second term becomes $u^K$, of the same
order as the first, and under the deep cancellation $|\tau|\ll|xy|+|c|$ the second term
dominates.
With $n=K{+}1$ the second term becomes a factor $u$ smaller and the condition under which
the result-relative reading dominates
\[
 |\tau|\;\gtrsim\;\bigl(Q^{\mathrm{full}}_K/P_K\bigr)\,u\,\bigl(|xy|+|c|\bigr),
\]
permits cancellation down to approximately $43.9$ bits for DW, $40.4$ bits for TW, and
$37.9$ bits for QW at $p=53$.

%\paragraph{Preservation of commutativity.}
The symmetrization stage (the leading \TS{} of the transposed pair) is independent of the
number of sweeps, and the expansion vector after symmetrization is position-wise
invariant under the transposition $x\leftrightarrow y$.
Because the added sweep is a deterministic map with a fixed order, the commutativity is
preserved in fully distilled form
(200{,}000 trials for each of the six types with no mismatch).

%\paragraph{When still deeper cancellation is required.}
Each additional sweep lowers the input-relative term by a further factor of $u$
($+6(m{-}1)$ flops), and deepens the admissible cancellation by approximately $p$ bits.
In the limit, it approaches the behavior of Priest's complete
non-overlapping normalization\cite{Priest91}.

%\paragraph{Numerical verification.}
Using exact rational arithmetic, we checked the identity (error $=$ the final truncation)
and the above bound (with the exact coefficients for $p=53$) on 8{,}000 trials each, of
which one-third was drawn from a deep-cancellation family (depth $5$--$45$ bits).
The identity held exactly in every trial, and there was no violation of the bound
(the maximum degree of saturation was $0.013$ for TW and $10^{-4}$ for QW).

%%% ===============================================================================
%%% Bibliography
%%% ===============================================================================

\end{document}

%% file: fig/dwnet_en.tex
\begin{tikzpicture}[
  gate/.style={draw,rounded corners,minimum width=15mm,minimum height=6mm,font=\scriptsize},
  op/.style={font=\scriptsize},
  hi/.style={-{Latex[length=1.5mm]}},
  lo/.style={-{Latex[length=1.5mm]},dashed},
  node distance=6mm and 10mm]
\node[gate] (tp) {$\TP(x_0,y_0)$};
\node[op,below left=2mm and 6mm of tp] (mul) {$P_{01}\!=\!x_0y_1,\ P_{10}\!=\!x_1y_0$};
\node[gate,right=16mm of tp] (ts1) {$\TS(P_{00},c_0)$};
\node[op,below=8mm of tp] (l) {$\ell\!=\!\mathrm{fl}(P_{01}\!+\!P_{10})$};
\node[op,right=4mm of l] (v) {$v\!=\!\mathrm{fl}(E_{00}\!+\!c_1)$};
\node[op,below=4mm of v] (w) {$w\!=\!\mathrm{fl}(v\!+\!\ell)$};
\node[op,below=4mm of ts1] (tp2) {$t_p\!=\!\mathrm{fl}(t\!+\!w)$};
\node[gate,right=12mm of tp2] (fts) {$\FTS(s,t_p)$};
\node[op,right=8mm of fts] (out) {$z=(z_0,z_1)$};
\draw[hi] (tp) -- node[op,above]{$P_{00}$} (ts1);
\draw[lo] (tp.south) to[bend right=10] node[op,left]{$E_{00}$} (v.north);
\draw[hi] (ts1) -- node[op,left]{$s$} (fts.north-|ts1) ;
\draw[lo] (ts1.south) -- node[op,right]{$t$} (tp2.north);
\draw[hi] (w) -- (tp2);
\draw[hi] (tp2) -- (fts);
\draw[hi] (ts1.south west) to[bend left=8] (fts);
\draw[hi] (fts) -- (out);
\end{tikzpicture}

%% file: fig/twnet_en.tex
\begin{tikzpicture}[
  lvl/.style={draw,rounded corners,minimum height=8mm,font=\scriptsize,align=left,
              fill=black!3,text width=62mm,inner sep=2mm},
  ren/.style={draw,rounded corners,font=\scriptsize,align=left,text width=42mm,inner sep=2mm},
  op/.style={font=\scriptsize},
  hi/.style={-{Latex[length=1.5mm]}},
  lo/.style={-{Latex[length=1.5mm]},dashed}]

\node[op] at (-0.9,0.95) {weight};
\node[op] at (-0.9,0.3)  {$u^0$};
\node[op] at (-0.9,-1.1) {$u^1$};
\node[op] at (-0.9,-2.6) {$u^2$};

\node[lvl,anchor=west] (l0) at (-0.4,0.3)  {$(B,r)\gets\TS(P_{00},c_0)$};
\node[lvl,anchor=west] (l1) at (-0.4,-1.1) {$A\gets\TS$ chain (3 gates): $(P_{01},P_{10})\to E_{00}\to c_1$\\
  $(m_1,m_2)\gets\TS(r,A)$\ \ {\tiny emits $q_{1..3}$}};
\node[lvl,anchor=west] (l2) at (-0.4,-2.6) {$G\gets$ plain sums:\\
  $\sigma=\mathrm{fl}(\mathrm{fl}(P_{02}+P_{20})+P_{11})$,\ \ $\mathrm{fl}(E_{01}+E_{10})$,\ \ $c_2$,\ \ $q_{1..3}$};

\node[ren,anchor=west] (rn) at (7.1,-1.1) {\textbf{Renormalize}\\
  $(w_0,w_1)\gets\FTS(B,m_1)$\\
  $(w_1,w_2)\gets\TS(w_1,m_2)$\\
  $(z_0,\rho)\gets\TS(w_0,w_1)$\\
  $(z_1,z_2)\gets\FTS(\rho,w_2)$};
\node[op,right=2mm of rn] (out) {$z$};

\draw[lo] (l0.south) -- node[op,right,pos=0.5]{$r$} (l1.north);
\draw[lo] (l1.south) -- node[op,right,pos=0.5]{$q_{1..3}$} (l2.north);
\draw[hi] (l0.east) to[out=0,in=160] node[op,above,pos=0.65]{$B$} (rn.north west);
\draw[hi] (l1.east) -- node[op,above]{$m_1$} (rn.west);
\draw[hi] (l2.east) to[out=0,in=200] node[op,below,pos=0.65]{$m_2$} (rn.south west);
\draw[hi] (rn) -- (out);
\end{tikzpicture}

%% file: fig/qwnet_en.tex
\begin{tikzpicture}[
  lvl/.style={draw,rounded corners,minimum height=8mm,font=\scriptsize,align=left,
              fill=black!3,text width=62mm,inner sep=2mm},
  ren/.style={draw,rounded corners,font=\scriptsize,align=left,text width=42mm,inner sep=2mm},
  op/.style={font=\scriptsize},
  hi/.style={-{Latex[length=1.5mm]}},
  lo/.style={-{Latex[length=1.5mm]},dashed}]

\node[op] at (-0.9,1.0)  {weight};
\node[op] at (-0.9,0.35) {$u^0$};
\node[op] at (-0.9,-1.1) {$u^1$};
\node[op] at (-0.9,-2.9) {$u^2$};
\node[op] at (-0.9,-4.6) {$u^3$};

\node[lvl,anchor=west] (l0) at (-0.4,0.35) {$(B,r)\gets\TS(P_{00},c_0)$};
\node[lvl,anchor=west] (l1) at (-0.4,-1.1) {$A_1\gets\TS$ chain (4 gates):\\
  $(P_{01},P_{10})\to E_{00}\to c_1\to r$\ \ {\tiny emits $f_{1..4}$}};
\node[lvl,anchor=west] (l2) at (-0.4,-2.9) {$A_2\gets\TS$ chain (9 gates):\\
  $(P_{02},P_{20})\to P_{11}\to \tilde{E}{=}\TS(E_{01},E_{10})$\\
  $\to c_2\to f_{1..4}$\ \ {\tiny emits $g_{1..9}$}};
\node[lvl,anchor=west] (l3) at (-0.4,-4.6) {$A_3\gets$ plain sums:\\
  $\mathrm{fl}(\mathrm{fl}(E_{02}{+}E_{20})+\mathrm{fl}(E_{11}{+}D))$,\ \ $c_3$,\ \ $g_{1..9}$};

\node[ren,anchor=west] (rn) at (7.1,-2.2) {\textbf{Renormalize}\\
  $(w_0,w_1)\gets\FTS(B,A_1)$\\
  $(w_1,w_2)\gets\TS(w_1,A_2)$\\
  $(w_2,w_3)\gets\TS(w_2,A_3)$\\
  $(z_0,\rho)\gets\TS(w_0,w_1)$\\
  $(z_1,\varsigma)\gets\TS(\rho,w_2)$\\
  $(z_2,z_3)\gets\FTS(\varsigma,w_3)$};
\node[op,right=2mm of rn] (out) {$z$};

\draw[lo] (l0.south) -- node[op,right,pos=0.5]{$r$} (l1.north);
\draw[lo] (l1.south) -- node[op,right,pos=0.5]{$f_{1..4}$} (l2.north);
\draw[lo] (l2.south) -- node[op,right,pos=0.5]{$g_{1..9}$} (l3.north);
\draw[hi] (l0.east) to[out=0,in=150] node[op,above,pos=0.6]{$B$} (rn.north west);
\draw[hi] (l1.east) to[out=0,in=175] node[op,above,pos=0.6]{$A_1$} (rn.west);
\draw[hi] (l2.east) to[out=0,in=195] node[op,below,pos=0.6]{$A_2$} (rn.south west);
\draw[hi] (l3.east) to[out=0,in=250] node[op,below,pos=0.7]{$A_3$} (rn.south);
\draw[hi] (rn) -- (out);
\end{tikzpicture}